\documentclass{article}

\usepackage[utf8]{inputenc}
\usepackage[english]{babel}

\usepackage[a4paper,left=.20\paperwidth,right=0.20\paperwidth]{geometry}

\usepackage{amsmath}
\usepackage{amsthm}
\usepackage{amssymb}
\usepackage{mathrsfs} 

\usepackage{xcolor}

\usepackage{graphicx}
\usepackage{subcaption} 

\usepackage{tikz}
\usetikzlibrary{arrows} 

\usepackage{array}
\usepackage{booktabs}
\usepackage{ragged2e}
\usepackage[font=small]{caption}

\usepackage{paralist}

\usepackage{accents}

\usepackage[linesnumbered, noend]{algorithm2e}
\newcommand{\algoAnd}{\\\quad{\bfseries \emph{and}}\xspace}
\let\oldnl\nl
\newcommand{\nonl}{\renewcommand{\nl}{\let\nl\oldnl}}
\newcommand{\var}{\textup}
\newcommand{\func}{\textsf}

\usepackage{xspace}

\usepackage{calc}

\usepackage{todonotes} 

\usepackage{zibtitlepage}

\usepackage{hyperref}

\input{commands}

\newcommand{\MIP}{\ensuremath{\mathscr{P}}\xspace}
\newcommand{\MIPs}{\ensuremath{\MIP^{\text{s}}}\xspace}
\newcommand{\MIPf}{\ensuremath{\MIP^{\text{f}}}\xspace}
\newcommand{\MIPsf}{\ensuremath{\MIP^{\text{sf}}}\xspace}

\newlength\myheight
\newlength\mydepth
\newcommand*\inlinegraphics[1]{%
  \settototalheight\myheight{Xygp}%
  \settodepth\mydepth{Xygp}%
  \raisebox{-\mydepth}{\includegraphics[height=\myheight]{#1}}%
}

\begin{document}

\newcommand{\myorcidlink}[1]{\,\href{https://orcid.org/#1}{\raisebox{-0.45ex}{\includegraphics[width=1.8ex]{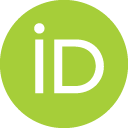}}}}

\ZTPAuthor{Felix~Hennings\protect\myorcidlink{0000-0001-6742-1983},
           Lovis~Anderson\protect\myorcidlink{0000-0002-4316-1862},
           Kai~Hoppmann-Baum\protect\myorcidlink{0000-0001-9184-8215},
           Mark~Turner\protect\myorcidlink{0000-0001-7270-1496},
           Thorsten~Koch\protect\myorcidlink{0000-0002-1967-0077}}
\ZTPTitle{Controlling transient gas flow in real-world pipeline intersection areas}
\ZTPInfo{The content of this report is also available as publication in Optimization and Engineering. Please always cite as:\\[1ex] Hennings, F., Anderson, L., Hoppmann-Baum, K., Turner, M., Koch, T. Controlling transient gas flow in real-world pipeline intersection areas. Optimization and Engineering, 2021, \href{https://doi.org/10.1007/s11081-020-09559-y}{\textcolor{blue}{DOI:10.1007/s11081-020-09559-y}}}
\ZTPNumber{19-24}
\ZTPMonth{June}
\ZTPYear{2019}

\title{Controlling transient gas flow in real-world pipeline intersection areas}
\author{Felix Hennings, Lovis Anderson, Kai Hoppmann-Baum, Mark Turner,\\
Thorsten Koch}
\zibtitlepage

\maketitle

\begin{abstract}
  Compressor stations are the heart of every high-pressure gas transport network.
  Located at intersection areas of the network they are contained in huge complex plants, where they are in combination with valves and regulators responsible for routing and pushing the gas through the network.
  Due to their complexity and lack of data compressor stations are usually dealt with in the scientific literature in a highly simplified and idealized manner.
  As part of an ongoing project with one of Germany's largest Transmission System Operators to develop a decision support system for their dispatching center, we investigated how to automatize control of compressor stations. Each station has to be in a particular configuration, leading in combination with the other nearby elements to a discrete set of up to 2000 possible feasible operation modes in the intersection area.
  Since the desired performance of the station changes over time, the configuration of the station has to adapt.
  Our goal is to minimize the necessary changes in the overall operation modes and related elements over time, while fulfilling a preset performance envelope or demand scenario.
  This article describes the chosen model and the implemented mixed integer programming based algorithms to tackle this challenge.
  By presenting extensive computational results on real world data we demonstrate the performance of our approach.
\end{abstract}

\section{Introduction}
  Throughout the past years, the mathematics of gas transport has been an intensively studied topic.
  While natural gas was, is, and will be one of the major energy sources in Germany, making the efficient and safe transport a field of high economical and political relevance~\cite{Fed2019}, the task is also challenging from a mathematical point of view.
  One such challenge can be found in the compressors, which push the gas through the network by increasing its pressure.
  Compressors are typically set up as a compressor station, whereby multiple compressor units can be placed in specific configurations and dynamically adjusted to meet the current needs, allowing for different compression ratios and flow rates.
  At intersections of major transportation pipelines, arrangements of multiple compressor stations as well as other elements like valves or pressure regulators can be found.
  Such an arrangement makes it possible to choose different connections of the intersecting pipelines, and operate the system for various flow directions and pressure levels.
  To optimize control of these areas taking all the technical restrictions into account is already combinatorially challenging.
  The complexity of the problem is further increased by the physics of gas flow.
  For pipes this physics is described by the Euler Equations~\cite{Osi1996}, a set of nonlinear hyperbolic partial differential equations (PDEs), which even in simplified versions yield computationally challenging constraints.

  Historically, research focused first on the simulation of gas flow, i.e., dealing with the partial differential equations given all the discrete decisions.
  This field has been studied for many decades already, see for example~\cite{BroGasHer2011} and the references therein.
  Over recent years, the optimization of gas transport including also the combinatorial aspects have gained more and more attention.
  In~\cite{RioBor2015} a general overview over optimization problems related to natural gas is given, which includes but is not restricted to the transport of gas.
  Most of the corresponding literature mentioned so far considers the stationary gas transport problem, which searches for one stable network state, making an algebraic description of the gas flow possible.
  An overview of state-of-the-art approaches for the stationary case can be found in~\cite{KocHilPfeSch2015} and~\cite{PfeFueGeiGei2014}, which consider large real-world instances and a huge amount of detail regarding the different network elements like compressors.

  This article deals with the more challenging variant of the problem: The transient gas transportation problem.
  Here, the goal is to find a set of control decisions on the elements over a future time horizon.
  For this problem, research is still in early stages.
  One of the first publications on transient gas transport optimization was~\cite{Mor2007}, who presented a mixed integer programming (MIP) model for the problem.
  In contrast to the structures described above, they only used a model consisting of single compressor units, whose compression capabilities are limited by a minimum and maximum power bound.
  This non-linear power bound as well as the non-linear pipe equations are approximated by piecewise linear functions.
  To solve the model for the objective of minimizing the compressor fuel costs, they used a special branching scheme for the piecewise linear functions as well as a simulated annealing heuristic described in~\cite{MahMarMor2007}.
  A little later,~\cite{DomGeiKolLan2011} also presented a solution to the problem of minimizing the compressor fuel costs.
  They modeled the problem according to~\cite{Mor2007}, using the same model for compressors and approximated the non-linear constraints by piecewise linear functions.
  However, they combined solving this MIP with solving a non-linear problem (NLP) formulation of the problem in an alternating way.
  From the solution of that non-linear problem they deduce a refinement of the piecewise linear approximations and repeat this procedure until finally arriving at a solution to the overall mixed-integer non-linear problem (MINLP) within a chosen approximation error.
  Other approaches tackle transient transport optimization problems, but neglect the discrete nature of some of the elements and therefore purely optimize over continuous variables, i.e., solve NLP problems.
  We mention as example the work of~\cite{ZloCheBac2015} and~\cite{MakVanZloHij2016}, who decide on the compression ratios of compressors, while again minimizing their fuel consumption.
  Very recently a few more studies on transient gas network optimization have been published.
  In~\cite{HahLeyZav2017} a specialized branching rule is used to solve a MINLP formulation of the problem with the objective to minimize fuel consumption.
  For the compressors they introduced the theoretical concept of different modes to switch between configurations, which each have separate feasible region.
  However, in the end they restrict to exactly one mode with nearly unrestricted compression capabilities for their experiments.
  Another approach combining different specialized solving techniques is presented in~\cite{GugLeuMarSch2018}, where iteratively a MIP model and a NLP model are solved for each single time step.
  These two models arise from the use of a special discretization of the Euler equations.
  For compression, both models use single compressor units featuring a linear feasible region.
  In contrast to the other mentioned publications, the objective function in~\cite{GugLeuMarSch2018} was not to minimize fuel cost, but to comply as well as possible with a set of future pressure and flow values given at the boundary nodes.
  Finally we mention~\cite{BurEggGroMar2019}, who considered maximizing the amount of temporarily stored gas in the network while maintaining a feasible transient control of the elements.
  They also introduce a new discretization of the Euler Equations, which results in a formulation close to the algebraic form of the stationary model.
  They then use this discretization to obtain globally optimal solutions.
  For the compressors however, again only single units restricted by upper and lower bounds in the compression ratio as well as the absolute pressure difference have been modeled.
  This problem is solved by alternating between solving a MIP model, which is obtained by replacing the non-linear constraints by piecewise linear functions, and a NLP model, in which the discrete decisions from the MIP are already fixed.

\begin{figure}[ht]
      \centering
      \includegraphics[width=\linewidth]{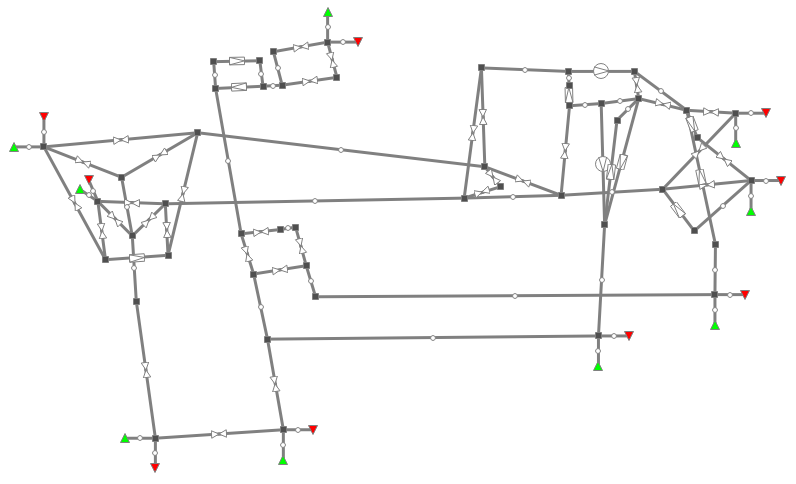}
      \caption[Example of the medium size network station E]{Example of the medium size network station E, see Table~\ref{tab:stationStatistics} for more details to its properties. The colored triangles represent the entry \inlinegraphics{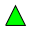} and exit \inlinegraphics{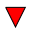} nodes of the station. Furthermore we have denoted the single network elements by \inlinegraphics{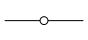} (pipe), \inlinegraphics{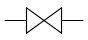} (valve), \inlinegraphics{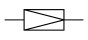} (regulator), and \inlinegraphics{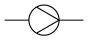} (compressor station).}
      \label{fig:gasNetworkStationExample}
  \end{figure}

  All approaches of the above mentioned publications on transient gas network optimization use a rather idealized model to deal with compressor stations.
  In contrast to this, we present in this article a transient gas network optimization problem featuring a compressor model with a so far unmatched amount of detail.
  The model is based on the modeling presented in~\cite{KocHilPfeSch2015} for the stationary case and takes the above described substructures into account.
  We focus on those network areas containing the compressor stations as well as additional active elements and call them \emph{network stations}, an example station is presented in Figure~\ref{fig:gasNetworkStationExample}.
  These areas contain the majority of active elements in the network.
  Regarding the number of contained elements they are comparable to, or even larger than the networks considered in the above mentioned literature on transient optimization.
  However, one difference to general gas network problems is the shortness of pipes due to the proximity of the elements in a network station, see Table~\ref{tab:stationStatistics} in Section~\ref{sec:computational_results} for an overview of different network stations.
  Because of their shortness, the pipes ability to store gas is negligible owing to their small volumes.
  Furthermore, the pressure loss induced by friction in the pipe is dependent on the its length and has therefore reduced impact.
  This allows us to use a linear pipe model as introduced in~\cite{Hen2018} without losing much accuracy and still producing realistic results.

  For each network station we are given an initial state as well as future demands in terms of both inflow and pressure levels at the boundaries.
  The goal is then to find a feasible control of all the network elements over time, which can be interpreted as a recommendation for network operators on how to control the network in the future.
  The overall objective is to meet the future demands as best as possible while simultaneously minimizing the total number of control changes.
  The latter is preferable since it reduces strain on the technical elements and enables the gas network operators to understand and actually perform the desired control recommendations.

  The rest of the article is structured as follows:
  In Section~\ref{sec:element_models} we will describe the mathematical models for all used elements and formulate a corresponding MIP model.
  The preprocessing needed to convert the given compressor data into a linear description of the feasible operating range will be introduced in Section~\ref{sec:compressor_configurations_definition}.
  However, solving the resulting MIP is quite challenging due to the complex compressor station model.
  We therefore propose a different solution approach in Section~\ref{sec:stationModelAlgorithm} based on solving slightly adjusted versions of the presented MIP. Section~\ref{sec:computational_results} then presents our results on computing solutions to a large number of real world networks situations.
  We finish with the conclusion in Section~\ref{sec:conclusion}.

\section{Mathematical model}
\label{sec:element_models}

  We model the gas network as a directed graph \generalGraph in which the arcs \setArcs represent the different network elements and nodes $\setVertices$ represent the junctions of the arcs.
  We split \setArcs into individual sets \setArcsAllStationModel for the network elements considered in this paper, i.e., pipes, valves, resistors, regulators, and compressor stations respectively.
  Note that regulators are also often named control valves in the literature, e.g. see~\cite{EhrSte2003} or~\cite{KocHilPfeSch2015}.
  In a similar fashion we split the node set \setVerticesAll into boundary nodes and inner nodes respectively.
  Here, boundary nodes \setBoundaryNodes represent those having inflow and pressure level demand values for future time steps.
  We define the set of considered time steps as \setTimestepsAll where \setTimestepsNoZeroAll are the future time steps having demand conditions at boundary nodes.
  Associated with each time step $t$ is a value \granularityI{t} representing the time difference in seconds from $t$ to the initial state time step 0.

  The most important quantities we will consider to describe the gas flow are pressure and mass flow.
  We have pressure variables \pressI{v,t} at each node $v\in\setVertices$ and time $t\in\setTimesteps$ as well as variables \mFlowI{a,t} for the flow from $l$ to $r$ on each non-pipe arc $(l,r)=a\in\setArcs\setminus\setPipes$ and time $t \in \setTimesteps$.
  For a pipe $(l,r)=a\in\setPipes$, we have two flow variables \mFlowI{l,a,t} representing the inflow into the pipe at end node $l$ and \mFlowI{r,a,t} representing the outflow out of the pipe at end node $r$, similar to the model of~\cite{EhrSte2003} on pipes or the one of~\cite{DomGeiKolLan2011} on all arcs.
  Our final flow quantity is the inflow at a boundary node \inflowI{v,t} entering the network for each $v\in\setBoundaryNodes$ and $t \in \setTimesteps$.
  Although we are given flow demands for the future, we allow deviations from them, and hence have to have a variable capturing the actual inflow value.

  We assume that there are bounds on all stated quantities for each point in time $t \in \setTimesteps$, so upper and lower bounds \pressUB{v,t} and \pressLB{v,t} on the pressure at each node $v\in\setVertices$, upper and lower bounds \mFlowUB{v,t} and \mFlowLB{v,t} on the flow on each non-pipe arc $a\in\setArcs\setminus\setPipes$ respectively the inflow and outflow of each pipe $a\in\setPipes$, as well as upper and lower bounds \paramInflowValueUB{v}{t} and \paramInflowValueLB{v}{t} on the inflow at each boundary node $v\in\setBoundaryNodes$.
  Note that while pressure is always positive, flow itself can be negative, as it can represent flow in the opposite direction. For example, negative inflow at a boundary node represents flow out of the network at this node.

  In the following we will describe each of the elements and at the same time introduce corresponding variables and constraints to create a MIP formulation.
  This MIP will not be solved directly, but is the basis for the three variants used in our overall solution algorithm described in Section~\ref{sec:stationModelAlgorithm}.

  Note that a list of all used variables can be found in the appendix in Table~\ref{tab:Variables}.

\subsection{Compressor stations}
\label{sec:compressor_station}

  Compressor stations are responsible for increasing the pressure in the network and thereby the most important in controlling the flow of gas.
  They are also the most complex elements, having their own substructure and a large amount of operational restrictions.
  Our model of the compressor station is based on the description in~\cite{FugGeiGolMor2015}, where the elements are called \emph{compressor groups} instead.

\paragraph{Structure of a compressor station}
  A compressor station $(l,r)=a\in\setCompressorStations$ has three different modes: Bypass, closed and active mode.
  In bypass mode the element is bypassed and therefore allows unrestricted gas flow without changing the pressure level.
  For the closed mode, the element is closed and thereby blocks the gas flow, which disconnects the network between the its end nodes.
  Finally, in active mode the gas is compressed and pressure is increased along the direction of the flow.

  When compressing, the compressor station can use a set of associated compressor units \setCompressorUnitsI{a}.
  These are the actual pressure increasing elements, each with a separate operating range.
  In the compressor station, these compressor units are combined in series and/or parallel to allow proper reactions to different compression requirements.
  The set of all allowed serial-parallel compressor unit combinations is called the set of configurations \setCompressorConfigurations{a} for a compressor station $a\in\setCompressorStations$, from which exactly one active configuration has to be chosen if the compressor station is in active mode.
  For each of these configurations $c\in\setCompressorConfigurations{a}$, we create a polytope in the space $(\pressI{l}, \pressI{r}, \mFlow)$ describing the feasible operating range of the compressor station using configuration $c$.
  This polytope is described as the intersection of a set of half spaces $\setConfigurationFacets{c} = \{(w,x,y,z) \in \mathbb{R}^4\}$ encoding inequalities of the form $w\cdot \pressI{l} + x\cdot \pressI{r} + y\cdot \mFlow + z \leq 0$.
  The creation of the feasible operating range of the configurations of the compressor stations is described in Section~\ref{sec:compressor_configurations_definition}.

\paragraph{Compressor station model}

  To model the above described constraints we use a disjunctive formulation.
  This is the most compact formulation in terms of number of constraints and variables for which holds that its LP relaxation is equal to the convex hull of its feasible points~\cite{Bal1985}\cite{Bal2018}.
  For this we introduce binary ``selection" variables \cfg{c,a,t} for each configuration $c\in\setCompressorConfigurations{a}$ of a compressor station $(l,r)=a\in\setCompressorStations$, as well as \modeBy{a,t} and \modeCl{a,t} for bypass and closed mode respectively.
  In addition we introduce corresponding sets of pressure and flow variables.
  The binary variables will force the pressure and flow variables of all non-selected configurations or modes to be zero and only enforce the constraints of the selected configuration or mode.
  The introduced pressure and flow variables are:
  \begin{align*}
    & \pressBy{a,t} \quad \mFlowBy{a,t} & \text{bypass mode variables} \\
    & \pressClL{a,t} \quad \pressClR{a,t} & \text{closed mode variables} \\
    & \pressCfgL{c,a,t}     \quad \pressCfgR{c,a,t} \quad \mFlowCfg{c,a,t} \quad \forall c\in\setCompressorConfigurations{a} & \text{configuration variables}
  \end{align*}
  Note that we need only one \press value for bypass mode, since here $\pressI{l} = \pressI{r}$ holds. Also there is no \mFlow variable for the closed mode, since $\mFlow=0$ holds in this case anyway.
  Furthermore, all introduced variables have bounds equal to the corresponding original pressure and flow bounds \pressUB{l,t}, \pressLB{l,t}, \pressUB{r,t}, \pressLB{r,t}, \mFlowUB{a,t}, \mFlowLB{a,t} of the compressor station $(l,r)=a\in\setCompressorStations$, enlarged to zero if necessary.
  We indicate these bounds by the variable symbol combined with an overscore respectively underscore.

  We are now able to state the constraints for all $(l,r) = a\in\setCompressorStations$ and $t\in\setTimestepsNoZero$:
  \newlength{\myLength}
  \setlength{\myLength}{2.5em}
  \begin{align}
    1 &= \sum_{c\in\setCompressorConfigurations{a}} \cfg{c,a,t} + \modeBy{a,t} + \modeCl{a,t} \label{eq:compressorStation_OneModeOrConfig} \\
    \pressI{l,t} &= \pressBy{a,t} + \pressClL{a,t} + \sum_{c\in\setCompressorConfigurations{a}} \pressCfgL{c,a,t} \label{eq:compressorStation_define_pL}\\
    \pressI{r,t} &= \pressBy{a,t} + \pressClR{a,t} + \sum_{c\in\setCompressorConfigurations{a}} \pressCfgR{c,a,t} \\
    \mFlowI{a,t} &= \mFlowBy{a,t} + \sum_{c\in\setCompressorConfigurations{a}} \mFlowCfg{c,a,t} \label{eq:compressorStation_last1}\\
    \pressCfgLLB{c,a,t}\cfg{c,a,t} \leq \makebox[\myLength][r]{\pressCfgL{c,a,t}} &\leq \pressCfgLUB{c,a,t}\cfg{c,a,t} \quad \forall c\in\setCompressorConfigurations{a} \label{eq:compressorStation_cfg_first}\\
    \pressCfgRLB{c,a,t}\cfg{c,a,t} \leq \makebox[\myLength][r]{\pressCfgR{c,a,t}} &\leq \pressCfgRUB{c,a,t}\cfg{c,a,t} \quad \forall c\in\setCompressorConfigurations{a} \\
    \mFlowCfgLB{c,a,t}\cfg{c,a,t} \leq \makebox[\myLength][r]{\mFlowCfg{c,a,t}} &\leq \mFlowCfgUB{c,a,t}\cfg{c,a,t} \quad \forall c\in\setCompressorConfigurations{a} \label{eq:compressorStation_cfg_last}\\
    \pressByLB{a,t}\modeBy{a,t} \leq \makebox[\myLength][r]{\pressBy{a,t}} &\leq \pressByUB{a,t}\modeBy{a,t} \label{eq:compressorStation_first2} \\
    \mFlowByLB{a,t}\modeBy{a,t} \leq \makebox[\myLength][r]{\mFlowBy{a,t}} &\leq \mFlowByUB{a,t}\modeBy{a,t} \\
    \pressClLLB{a,t}\modeCl{a,t} \leq \makebox[\myLength][r]{\pressClL{a,t}} &\leq \pressClLUB{a,t}\modeCl{a,t} \\
    \pressClRLB{a,t}\modeCl{a,t} \leq \makebox[\myLength][r]{\pressClR{a,t}} &\leq \pressClRUB{a,t}\modeCl{a,t} \label{eq:compressorStation_last2}\\
    w\cdot \pressCfgL{c,a,t} + x\cdot \pressCfgR{c,a,t} + y\cdot \mFlowCfg{c,a,t} + z\cfg{c,a,t} &\leq 0 \quad \forall (w,x,y,z) \in\setConfigurationFacets{c} \quad \forall c\in\setCompressorConfigurations{a} \label{eq:compressorStation_cfg_facets}
  \end{align}

\subsection{Pipes}
\label{sec:pipes}

  Gas flow in pipelines is for operational purposes modeled as one-dimensional flow through a straight cylindrical pipe.
  When assuming a constant gas temperature $T$,
  the flow can be described by the isothermal \emph{Euler Equations}\cite{Osi1996} consisting of the \emph{Continuity Equation} and the \emph{Momentum Equation}.
  The Continuity Equation describes the conservation of mass, i.e., guaranteeing that mass may neither be created nor destroyed.
  On the other hand, the Momentum Equation reflects the equality between the force acting on the gas particles and their corresponding rate of change of momentum.
  For a pipe $a=(l,r)$ we can state this pair of equations as
  \begin{align}
    \frac{\partial \dens}{\partial t} + \frac{\partial (\dens\velo)}{\partial x} &=0 \label{eq:euler_continuity}\\
    \frac{\partial (\dens\velo)}{\partial t} + \frac{\partial \press}{\partial x} + \frac{\partial (\dens\velo^2)}{\partial x} + \frac{\fricI{a}}{2\diamI{a}}|\velo|\velo\dens + \gravAcc\slopeI{a}\dens &= 0 \label{eq:euler_momentum}.
  \end{align}
  Here $x$ denotes the position in the pipe by its distance from the source node $l$,
  $t$ the current time, \dens the density of the gas, \velo its velocity, \diamI{a} the diameter of the pipe, \gravAcc the gravitational acceleration, and \fricI{a} the friction factor of the pipe, which we assume to depend on pipe characteristics only, see Section~\ref{sec:frictionAndCompressibility}.
  With $\slopeI{a}\in\left[-1,1\right]$ we denote the slope of the pipe, i.e., the quotient of the elevation increase between the pipes endpoints and the length \lenI{a} of the pipe.

  In order to complete the system of equations describing the state variables \press, \dens, and \velo, we add the \emph{equation of state for real gases} to establish the connection between \press and \dens as
  \begin{equation}
    \press = \dens \sGasConst \temp \zFactI{a}.
    \label{eq:realGasState}
  \end{equation}
  The two new quantities that arise are the compressibility factor \zFactI{a} and the specific gas constant  \sGasConst.
  We assume both values to be constant parameters, which for the compressibility factor is a common assumption in the gas transport literature, see for example~\cite{Osi1996}\cite{BurEggGroMar2019}.
  For the specific gas constant this follows from its dependence on the molar mass, which in turn is determined by the gas mixture which we assume to be constant.

  In the following, we will drop the terms $\partial_t(\dens\velo)$ and $\partial_x (\dens\velo^2)$ as they contribute only little to the equation under normal operating conditions\,\cite{EhrSte2003}\cite{Osi1996}.
  In addition, we reformulate the pipe flow equations in terms of the quantities we are interested in, i.e., pressure \press and mass flow \mFlow, where \mFlow is defined using the cross sectional area $\areaI{a}=\diamI{a}^2\frac{\pi}{4}$ of the cylindric pipe $a$ as
  \begin{equation}
    \mFlow = \areaI{a}\dens\velo. \label{eq:massflow_definition}
  \end{equation}
  Then we can write \eqref{eq:euler_continuity} and \eqref{eq:euler_momentum} as
  \begin{align*}
    \frac{\partial \press}{\partial t}
    + \frac{\sGasConst \temp \zFactI{a}}{\areaI{a}}\frac{\partial \mFlow}{\partial x}
    &=0 \\
      \frac{\partial \press}{\partial x}
    + \frac{\fricI{a}\sGasConst\temp\zFactI{a}}{2\diamI{a}\areaI{a}^2}\frac{|\mFlow|\mFlow}{\press}
    + \frac{\gravAcc\slopeI{a}}{\sGasConst\temp\zFactI{a}}\press
    &= 0.
  \end{align*}
  Since the spatial pressure change now mainly depends on the friction term including the friction factor $\fricI{a}$ this model variant is often referred to as the \emph{friction dominated} model, see model (FD1) in~\cite{BroGasHer2011} respectively model (ISO3) in~\cite{DomHilLanTis2017}.

  For the discretization, we use the implicit box scheme introduced by~\cite{DomGeiKolLan2011}, respectively~\cite{KolLanBal2010}.
  Here, the spacial domain is the length $\lenI{a}$ of pipe $a=(l,r)$ and the time domain is the set of time steps \setTimesteps, which we defined in the beginning of Section~\ref{sec:element_models}.
  Using the notation of flow into and out of a pipe, again defined in beginning of Section~\ref{sec:element_models}, as well as the function \granularity, we are able to write the discretized model for two adjacent time points $t_0$ and $t_1$ as
  \begin{align}
    \pressI{l,t_1} + \pressI{r,t_1} - \pressI{l,t_0} - \pressI{r,t_0}
    + \frac{2\sGasConst\temp\zFactI{a}(\granularityI{t_1} - \granularityI{t_0})}{\lenI{a}\areaI{a}}
      \left(\mFlowI{r,a,t_1} - \mFlowI{l,a,t_1}\right) &=0 \label{eq:pipes_dicretized_continuity}\\
    \pressI{r,t_1} - \pressI{l,t_1}
    + \frac{\fricI{a}\sGasConst\temp\zFactI{a}\lenI{a}}{4\diamI{a}\areaI{a}^2}
      \left(  \frac{|\mFlowI{l,a,t_1}|\mFlowI{l,a,t_1}}{\pressI{l,t_1}}
            + \frac{|\mFlowI{r,a,t_1}|\mFlowI{r,a,t_1}}{\pressI{r,t_1}}\right)& \notag\\
    + \frac{\gravAcc\slopeI{a}\lenI{a}}{2\sGasConst\temp\zFactI{a}}
      \left( \pressI{l,t_1} + \pressI{r,t_1} \right) &=0. \label{eq:pipes_dicretized_momentum}
  \end{align}

  As a final step, we will linearize the Momentum Equation~\eqref{eq:pipes_dicretized_momentum} as proposed in~\cite{Hen2018} by fixing the absolute velocity $|\velo|$ to a predefined constant in the friction term, where the absolute velocity is according to \eqref{eq:realGasState} and \eqref{eq:massflow_definition} defined as
  $|\velo|=|\frac{\sGasConst\temp\zFactI{a}}{\areaI{a}}\frac{\mFlow}{\press}|
  =\frac{\sGasConst\temp\zFactI{a}}{\areaI{a}}\frac{|\mFlow|}{\press}$.
  We found that the relative error stated in~\cite{Hen2018} is less relevant for the overall accuracy of the model, if the pipes are contained in network stations.
  The reason for this is, that the elements of a network station are usually clustered in a small geographic area.
  Therefore the contained pipes are relatively short in comparison to the rest of the network, see also Table~\ref{tab:stationStatistics} with statistics over different network stations.
  Since the friction based pressure reduction depends on the pipe length, the corresponding term and therefore the also stated relative error have much less impact than usual.
  This linearization allows us to model the pipe flow in a MIP context.
  The final equations for each pipe $(l,r)=a\in\setPipes$ and all adjacent time points $t_0,t_1\in\setTimesteps$ are:
  \begin{align}
    \pressI{l,t_1} + \pressI{r,t_1} - \pressI{l,t_0} - \pressI{r,t_0}
    + \frac{2\sGasConst\temp\zFactI{a}(\granularityI{t_1} - \granularityI{t_0})}{\lenI{a}\areaI{a}}
      \left(\mFlowI{r,a,t_1} - \mFlowI{l,a,t_1}\right) &=0 \label{eq:pipes_constVelo_continuity}\\
    \pressI{r,t_1} - \pressI{l,t_1}
  + \frac{\fricI{a}\lenI{a}}{4\diamI{a}\areaI{a}}
    \left(\absVelo{l,a}\mFlowI{l,a,t_1} + \absVelo{r,a}\mFlowI{r,a,t_1}\right) & \notag\\
  + \frac{\gravAcc\slopeI{a}\lenI{a}}{2\sGasConst\temp\zFactI{a}}
    \left(\pressI{l,t_1} + \pressI{r,t_1}\right) &=0 \label{eq:pipes_constVelo_momentum}
  \end{align}
  The constant \absVelo{x,a} for one of the end nodes $x\in\{l,r\}$ of pipe $a=(l,r)$ is determined based on the flow and pressure values of the given initial state, i.e.,
  \begin{equation*}
    \absVelo{x,a} = \frac{\sGasConst\temp\zFactI{a}}{\areaI{a}}\frac{|\mFlowI{x,a,0}|}{\pressI{x,0}}.
  \end{equation*}

  \subsubsection{Friction and compressibility factor}
  \label{sec:frictionAndCompressibility}
  The \emph{Darcy–Weisbach friction factor} \fric, describes the pressure drop on a pipe $a$ caused by frictional forces and depends on the diameter \diamI{a} and integral roughness \roughI{a} of the pipe, as well as the current flow \mFlow and the dynamic viscosity \dynVisc of the gas.
  For turbulent gas flow the most accurate description is given by the implicit Colebrook-White equation~\cite{ColWhi1937}\cite{FugGeiGolMor2015}.
  There exist a series of different explicit approximation formulas, typically depending on the \emph{Reynolds Number} which describes the amount of turbulence of the flow.
  We use the formula of Nikuradse~\cite{Nik1950}\cite{FugGeiGolMor2015}, which assumes infinite turbulence and makes the friction factor dependent only on the constant diameter \diamI{a} and integral roughness \roughI{a} of the pipe:
  \begin{equation*}
    \fricI{a} = \left( 2 \log_{10} \left(\frac{\diamI{a}}{\roughI{a}}\right)+1.138\right)^{-2}.
  \end{equation*}

  For the compressibility factor $\zFact$ we use the approximation formula developed by Papay~\cite{Pap1968}\cite{Sal2002}, which is valid up to 150\,bar and thus fits our considered pressure range of 1\,bar to 100\,bar well.
  It is given as
  \begin{equation*}
    \zFact(\press) = 1 - 3.52 \frac{\press}{\pseudoP} e^{-2.26 \frac{\temp}{\pseudoT}} + 0.247 \left(\frac{\press}{\pseudoP}\right)^2 e^{-1.878 \frac{\temp}{\pseudoT}}.
  \end{equation*}
  Apart from the pressure \press and gas temperature \temp Papay's formula depends on the gas mixture dependent pseudo-critical pressure \pseudoP and temperature \pseudoT, which we assume to be given constants.
  The constant compressibility factor \zFactI{a} per pipe $a=(l,r)$ is then determined as an average of the corresponding values derived from the initial state pressure values of the pipes end nodes, i.e., by
  $\zFactI{a}= (\zFact(\pressI{l,0}) + \zFact(\pressI{r,0})) / 2$.

\subsection{Resistors}
\label{sec:resistors}

  Resistors are artificial elements created to model points of high friction in the network, which can be caused by all sorts of special elements.
  Some examples of these elements are measuring equipment and complex local piping, which both are not captured by the other considered element types but need to be accounted for.
  The pressure drop caused by a resistor arc $(l,r) = a \in\setResistors$ for time $t\in\setTimestepsNoZero$ is defined by the \emph{Darcy-Weisbach equation}~\cite{FugGeiGolMor2015}:
  \begin{equation*}
      \pressI{l,t} - \pressI{r,t}
    = \frac{\drag{a}\sGasConst\temp\zFactI{a}}{2\areaI{a}^2}
      \left( \frac{|\mFlowI{a,t}|\mFlowI{a,t}}{\pressI{\text{in},t}} \right)
  \end{equation*}
  Here the friction factor of the resistor is called the \emph{drag factor}, and is represented by \drag{a}, a parameter of the element.
  The compressibility factor \zFactI{a} is determined in the same way as described for pipes, see Section~\ref{sec:frictionAndCompressibility}.
  The formula is flow direction dependent, where \pressI{\text{in},t} is either \pressI{l,t} or \pressI{r,t} depending on \mFlowI{a,t} being positive or negative (for $\mFlowI{a,t} = 0$ it holds that $\pressI{l,t} = \pressI{r,t}$).

  As we did for pipes in Equation~\eqref{eq:pipes_constVelo_momentum} we linearize the model by assuming a constant velocity
  $|\velo|=\frac{\sGasConst\temp\zFactI{a}}{\areaI{a}}\frac{|\mFlow|}{\press}$, which also includes the flow direction dependent pressure value. The equations for each arc $(l,r) = a \in\setResistors$ and time $t\in\setTimestepsNoZero$ then reads as
  \begin{equation}
      \pressI{l,t} - \pressI{r,t}
    = \frac{\drag{a}\absVelo{a}}{2\areaI{a}}\mFlowI{a,t}. \label{eq:resistors_constVelo}
  \end{equation}
  The constant velocity value is again calculated based on the initial element state, and is defined as an average of the two velocities using the pressure from the corresponding resistors end nodes as
  \begin{align*}
    \absVelo{l} &= \frac{\sGasConst\temp\zFactI{a}}{\areaI{a}}\frac{|\mFlowI{a,0}|}{\pressI{l,0}} \\
    \absVelo{r} &= \frac{\sGasConst\temp\zFactI{a}}{\areaI{a}}\frac{|\mFlowI{a,0}|}{\pressI{r,0}} \\
    \absVelo{a} &:= \frac{\absVelo{l} + \absVelo{r}}{2}
  \end{align*}

\subsection{Valves}

  Valves are active elements that dynamically connect or disconnect two nodes by being open or closed respectively, and thereby change the network topology.
  Closed valves work exactly as the closed mode of a compressor station, while open valves imply the same behavior as the corresponding stations bypass mode.
  The mode of a valve is captured by the binary variable \modeOp{a,t}.
  For a valve arc $(l,r) = a \in\setValves$ and time $t\in\setTimestepsNoZero$ we write
  \begin{align}
    \pressI{l,t} - \pressI{r,t} &\leq ( 1 - \modeOp{a,t} )(\pressUB{l,t}-\pressLB{r,t}) \label{eq:valves_first}\\
    \pressI{l,t} - \pressI{r,t} &\geq ( 1 - \modeOp{a,t} )(\pressLB{l,t}-\pressUB{r,t}) \\
    \mFlowI{a,t} &\leq ( \modeOp{a,t} )\mFlowUB{a,t} \\
    \mFlowI{a,t} &\geq ( \modeOp{a,t} )\mFlowLB{a,t}. \label{eq:valves_last}
  \end{align}

\subsection{Regulators}

  A regulator or control valve is a valve with variable opening degree, used to reduce the pressure along the direction of the flow.
  Regulators can reduce the pressure in active mode and also be bypassed or closed like a compressor station.
  For each mode there is a binary variable, from which exactly one is equal to 1 at any time, i.e., the regulator always has to have a unique mode
  \begin{align}
    1 &= \modeCl{a,t} + \modeBy{a,t} + \modeAc{a,t} \qquad \forall a \in \setRegulators \quad \forall t \in \setTimestepsNoZero \label{eq:regulator_mode}
  \end{align}
  The implications of each of the modes can be modeled by the following constraints for each arc $a \in \setRegulators$ and all times $t \in \setTimestepsNoZero$
  \begin{align}
    \pressI{l,t} - \pressI{r,t} &\leq  + ( 1 - \modeBy{a,t} )(\pressUB{l,t}-\pressLB{r,t}) \label{eq:regulator_mode_definition_first} \\
    \pressI{l,t} - \pressI{r,t} &\geq  + ( 1 - \modeBy{a,t} - \modeAc{a,t} )(\pressLB{l,t}-\pressUB{r,t}) \\
    q_a &\leq ( 1 - \modeCl{a,t} )\mFlowUB{a,t} \\
    q_a &\geq 0. \label{eq:regulator_mode_definition_last}
  \end{align}
  Note that the flow is always positive, even in bypass mode. This occurs because all regulators have a flap trap which prevents flow going against the topological orientation.

\subsection{Nodes}

  The nodes don't represent technical elements but rather establish the connections between them.
  The pressure coupling is realized by using the pressures at a node in the constraints of all its incident arcs.
  To connect the mass flow between arcs we have flow conservation constraints in each node.
  This means that the sum of incoming flows should match the sum of outgoing flows, resulting in the following constraints for all $t\in\setTimestepsNoZero$:
  \begin{align}
     & \sum_{(l,v)=a\in\setPipes} \mFlowI{v,a,t} - \sum_{(v,r)=a\in\setPipes} \mFlowI{v,a,t} \notag\\
    +& \sum_{(l,v)=a\in\setArcs\setminus\setPipes} \mFlowI{a,t} - \sum_{(v,r)=a\in\setArcs\setminus\setPipes} \mFlowI{a,t} + \inflowI{v,t} = 0 \qquad \forall v\in\setBoundaryNodes \label{eq:flowConservation_boundaryNodes}\\
     & \sum_{(l,v)=a\in\setPipes} \mFlowI{v,a,t} - \sum_{(v,r)=a\in\setPipes} \mFlowI{v,a,t} \notag\\
    +& \sum_{(l,v)=a\in\setArcs\setminus\setPipes} \mFlowI{a,t} - \sum_{(v,r)=a\in\setArcs\setminus\setPipes} \mFlowI{a,t} = 0 \qquad \forall v\in\setInnerNodes \label{eq:flowConservation_innerNodes}
  \end{align}

\subsection{Network station}
\label{sec:gasNetworkStationDefinition}

  In addition to the constraints imposed by the single elements used in the network, the network station itself has a set of associated constraints.
  The model described below is an extension of the one described in~\cite{FugGeiGolMor2015} for similar structures.
  There, network stations are called \emph{compressor stations} and the elements we call compressor stations are called \emph{compressor groups}, as already mentioned in Section~\ref{sec:compressor_station}.

\paragraph{Network station structure}
  Most important are the \emph{operation modes} \setNSModes of the network station from which exactly one has to be selected at each time point and which determine the modes and configurations of all valves and compressor stations.
  Note that not all possible mode combinations of the different elements have to be valid operation modes of the network station.
  In addition to the operation mode a \emph{flow direction} has to be chosen for the network station from the set of possible flow directions \setFlowDirectionsNavi.
  The restrictions on each flow direction of the station prescribe the flow patterns in terms of inflow, outflow or no-flow over the boundary nodes of the station.
  The flow direction itself also has to fit to the selected operation mode of the network station, where the set of feasible operation mode and flow direction pairs is \setNSValidPairs.

\paragraph{Operation modes model}
  We introduce binary variables \opMode{o,t} for each $o\in\setNSModes$ and $t \in \setTimesteps$. They represent whether the operation mode has been selected, or respectively if the network station is in the given operation mode at this point in time.
  Furthermore, we define the function $M(o,a)$, which maps operation modes to the individual modes and configurations of valves and compressor stations:
  \begin{align*}
    M(o,a) := x &\text{ where }x\text{ is the mode or configuration of arc }a\\
    & \text{ in operation mode }o \quad \forall o\in\setNSModes\quad\forall a\in\setValves\cup\setCompressorStations\\
    \text{with}\qquad & x\in \{\text{op}, \text{cl}\} \text{ if }a\in\setValves \\
                      & x\in \{\text{by}, \text{cl}\} \cup \setCompressorConfigurations{a} \text{ if }a\in\setCompressorStations
  \end{align*}
  Note that we assume that all valves and compressor stations are determined by the network station operation modes.
  In general however, there are valves which are not controlled by the operation modes, and whose mode is already given as a fixed decision over time and cannot be changed.
  We can handle these valves by pre-processing, and either contracting the associated valve in the open case, or simply removing the valve in the closed case.

  Using $M(o,a)$ we can then state the operation mode related constraints for all $t\in\setTimestepsNoZero$:
  \begin{align}
    \sum_{o \in \setNSModes} \opMode{o,t} &= 1 \label{eq:opMode_choice} \\
    \modeOp{a,t}   &= \sum_{o\in\setNSModes : M(o,a)=\text{op}} \opMode{o,t} \quad\forall a\in\setValves \label{eq:valve_opMode_coupling}\\
    \modeBy{a,t} &= \sum_{o\in\setNSModes : M(o,a)=\text{by}} \opMode{o,t} \quad\forall a\in\setCompressorStations \label{eq:cs_bypass_opMode_coupling} \\
    \modeCl{a,t} &= \sum_{o\in\setNSModes : M(o,a)=\text{cl}} \opMode{o,t} \quad\forall a\in\setCompressorStations \label{eq:cs_closed_opMode_coupling} \\
    \cfg{c,a,t} &= \sum_{o\in\setNSModes : M(o,a)=c}   \opMode{o,t} \quad\forall c\in\setCompressorConfigurations{a}\quad\forall a\in\setCompressorStations \label{eq:cs_cfg_opMode_coupling}
  \end{align}
  Note that either Equation~\eqref{eq:cs_bypass_opMode_coupling}, Equation~\eqref{eq:cs_closed_opMode_coupling} or one of the constraints of Equation~\eqref{eq:cs_cfg_opMode_coupling} for one of the configurations $c\in\setCompressorConfigurations{a}$ can be omitted, because it follows from the remaining constraints combined with Equation~\eqref{eq:compressorStation_OneModeOrConfig}.

\paragraph{Operation mode unavailability}
  Certain operation modes are not available at specific points in time.
  The basis for this is the non-availability of compressor units over time, which is part of the model input data.

  As explained in Section~\ref{sec:compressor_station}, a configuration $c\in\setCompressorConfigurations{a}$ of some compressor station $a\in\setCompressorStations$ represents the serial and/or parallel combination of a subset of the compressor stations compressor units.
  Hence, the unavailability of a certain compressor unit at time $t$ results in the unavailability of all configurations which use this unit.
  On the next level, each network station operation mode defines the mode and (for the active mode) the configuration of each compressor station in the network station.
  Hence, all network station operation modes using a configuration for a compressor station which is unavailable for time $t$ will be unavailable for $t$, too.
  To implement this in the model, we just fix the variables \opMode{o,t} for the corresponding operation mode $o$ and time points $t$ to zero, i.e., remove them from the model.

  The unavailability of a compressor unit may not be aligned with the set of discrete time points \setTimesteps, i.e., the unavailability period may start or stop in between two adjacent time points.
  To be able to tell which of the two points is then effected by this, we have to establish an interpretation for the operation mode of a network station between two time steps.
  Therefore, we define that if a network station has the operation mode $A$ at the discrete time point $t$, then we also assume the station to have operation mode $A$ in the following time interval up to the next discrete time point $t+1$.
  From this definition it follows that if a network station operation mode is unavailable for some $k \in (t, t+1)$ with $t, t+1\in\setTimesteps$, then the operation mode is unavailable for time $t$ but potentially available for time $t+1$.

\paragraph{Operation mode transition times}
  \newcommand{\transTime}[2]{\ensuremath{\theta(#1,#2)}\xspace}

  If the operation mode of a network station is changed from mode $A$ to some other mode $B$, the transition takes a given amount of time \transTime{A}{B}. A transition time is given for every possible combination of operation modes and is part of the input data.
  While in transition between the two modes, the network station acts as follows:
  Assume the transition starts at time $t_0$, then for $t\in\left[t_0,t_0 + \frac{\transTime{A}{B}}{2}\right)$ the station uses mode $A$ while for $t\in\left[t_0 + \frac{\transTime{A}{B})}{2},t_0 + \transTime{A}{B}\right)$ the station uses mode $B$.
  In other words, the station stays in mode $A$ until it reaches the middle of the transition period and then changes to mode $B$.
  While being in transition, the whole transition period is blocked for other changes, i.e., two transition periods cannot overlap.
  Since we are only able to change network station operation modes at discrete time points, we assume that for each transition the middle point is in \setTimesteps.
  This is also in line with our interpretation of network station operation modes in between discrete time points.

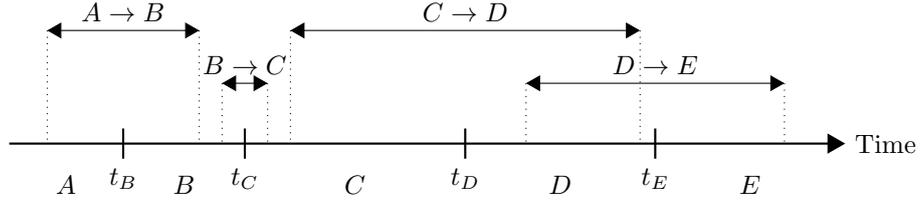
\begin{figure}[ht]
  \begin{tikzpicture}[>=triangle 60]
    \newcommand\XZ{1.5}   
    \newcommand\XI{3.1}   
    \newcommand\XII{6}    
    \newcommand\XIII{8.5} 
    \newcommand\XEND{11}
    \newcommand\TAB{1.0}
    \newcommand\TBC{0.3}
    \newcommand\TCD{2.3}
    \newcommand\TDE{1.7}
    \draw[thick,black,->] (0,0)  -- (\XEND,0) node[right] {Time};
    \draw[thick,black] (\XZ,  0.2)  -- (\XZ,  -0.2) node[below] {$t_B$};
    \draw[thick,black] (\XI,  0.2)  -- (\XI,  -0.2) node[below] {$t_C$};
    \draw[thick,black] (\XII, 0.2)  -- (\XII, -0.2) node[below] {$t_D$};
    \draw[thick,black] (\XIII,0.2)  -- (\XIII,-0.2) node[below] {$t_E$};
    \node[below] at (\XZ/2,            -0.3) {$A$};
    \node[below] at (\XZ/2   + \XI/2,  -0.3) {$B$};
    \node[below] at (\XI/2   + \XII/2, -0.3) {$C$};
    \node[below] at (\XII/2  + \XIII/2,-0.3) {$D$};
    \node[below] at (\XIII/2 + \XEND/2, -0.3) {$E$};
    \draw[<->] (\XZ   - \TAB,1.5)  -- (\XZ   + \TAB,1.5) node[midway,above] {$A\rightarrow B$};
    \draw[<->] (\XI   - \TBC,0.8)  -- (\XI   + \TBC,0.8) node[midway,above] {$B\rightarrow C$};
    \draw[<->] (\XII  - \TCD,1.5)  -- (\XII  + \TCD,1.5) node[midway,above] {$C\rightarrow D$};
    \draw[<->] (\XIII - \TDE,0.8)  -- (\XIII + \TDE,0.8) node[midway,above] {$D\rightarrow E$};
    \draw[dotted] (\XZ   - \TAB,1.5)  -- (\XZ   - \TAB,0.0);
    \draw[dotted] (\XZ   + \TAB,1.5)  -- (\XZ   + \TAB,0.0);
    \draw[dotted] (\XI   - \TBC,0.8)  -- (\XI   - \TBC,0.0);
    \draw[dotted] (\XI   + \TBC,0.8)  -- (\XI   + \TBC,0.0);
    \draw[dotted] (\XII  - \TCD,1.5)  -- (\XII  - \TCD,0.0);
    \draw[dotted] (\XII  + \TCD,1.5)  -- (\XII  + \TCD,0.0);
    \draw[dotted] (\XIII - \TDE,0.8)  -- (\XIII - \TDE,0.0);
    \draw[dotted] (\XIII + \TDE,0.8)  -- (\XIII + \TDE,0.0);
  \end{tikzpicture}
  \caption{Transition time example with 5 operation modes and 4 transitions, from which two are in conflict. The time point $t_X$ represent the first discrete time point in \setTimesteps in which mode $X$ is active.}
  \label{fig:transitionTimes}
  \end{figure}

  In Figure~\ref{fig:transitionTimes} we see an example of an operation mode sequence with corresponding transition times.
  The time point $t_X$ represents the first $t\in\setTimesteps$ in which operation mode $X$ is active, which according to the interpretation in the previous paragraph is also the first discrete time point $t$ in which $X$ is active.
  In this example, there would be a conflict in the transition times \transTime{C}{D} from mode $C$ to mode $D$ and \transTime{D}{E} from mode $D$ to mode $E$, since the two time periods overlap.
  Note that this is true, although for each single transition period the network station has the correct operation mode for each point in time, i.e. first mode in the first half of the period and the second mode in the second half.

  We will \emph{not} cover the transition time restrictions in the MIP model, but will make sure our solutions respect them in a different way, see Section~\ref{sec:determine_operation_modes}.
  For this we only need a way to check for a given sequence of operation modes if all the transition times are valid, i.e., the corresponding periods do not overlap, which we do as follows:
  For each network station operation mode in the sequence, we check if the time in which the operation mode is active is at least as big as the sum of adjacent transition time components, i.e., the sum of half the time of the transition into that mode and half of the time of the transition out of that mode.
  In the given example, we would have to check for mode $C$ if $\granularityI{t_D} - \granularityI{t_C} \geq \frac{\transTime{B}{C}}{2} + \frac{\transTime{C}{D}}{2}$ holds, where \granularityI{t} represents the time difference of a time step $t$ from the initial state time.
  If the above holds for all operation modes in the sequence, then the transition periods can never overlap since they are centered around the time points in which the mode change happens.

  Note, that for the last mode in the sequence, we do not need to do any checks at all.
  We do so as we cannot determine how long into the future the operation mode will be active for.
  Hence, we assume it is active long enough to comply with the given transition time, allowing the possibility to apply operation mode changes up to the very last time step.
  For the first mode, a similar assumption would be too optimistic, since the desired operation mode change would have then already been triggered before the start our time horizon.
  However, there is no known transition into this first mode, so we need only check if half the transition time into the next mode fits into the first modes active period.

\paragraph{Flow directions model}
  Similar to the way we handle operations modes, we introduce binary variables \flowDir{\flowDirection,t} representing the selection of flow direction $\flowDirection\in\setFlowDirectionsNavi$ at time $t \in \setTimesteps$.
  To be able to state the connection between the chosen flow direction and the actual boundary node inflow pattern, we represent each flow direction $\flowDirection$ as a tuple of the set of boundary nodes having inflow into the station $\flowDirectionEntries$ and the set of boundary nodes having outflow out of the station $\flowDirectionExits$.
  Hence, using the power set \powerset, a flow direction $\flowDirection$ is defined as
  \begin{equation*}
    (\flowDirectionEntries, \flowDirectionExits) = \flowDirection \in \setFlowDirectionsNavi \subseteq \powerset(\setBoundaryNodes) \times \powerset(\setBoundaryNodes) \text{ where } \flowDirectionEntries \cap \flowDirectionExits = \emptyset
  \end{equation*}
  Note that if $v\not\in\flowDirectionEntries$ and $v\not\in\flowDirectionExits$ for flow direction $(\flowDirectionEntries, \flowDirectionExits)$, then the inflow of node $v$ is zero.

  Using the inflow variable \inflowI{v,t} for boundary node $v$ and time $t$ we can define the flow direction constraints for each $t\in\setTimestepsNoZero$ as:
  \begin{align}
    \sum_{\flowDirection \in \setFlowDirectionsNavi} \flowDir{f,t} &= 1 \label{eq:flowDir_choice}\\
    \opMode{o,t} &\leq \sum_{(o,\flowDirection)\in\setNSValidPairs} \flowDir{\flowDirection,t} \quad \forall o\in\setNSModes \label{eq:ns_om_fd_coupling}\\
    \inflowI{v,t} &\geq ( 1 - \sum_{\flowDirection=(\flowDirectionEntries, \flowDirectionExits)\in\setFlowDirectionsNavi : v\not\in\flowDirectionExits}  \flowDir{\flowDirection,t})\paramInflowValueLB{v}{t} \quad\forall v\in\setBoundaryNodes \label{eq:ns_fenceNodes_flowDirection1}\\
    \inflowI{v,t} &\leq ( 1 - \sum_{\flowDirection=(\flowDirectionEntries, \flowDirectionExits)\in\setFlowDirectionsNavi : v\not\in\flowDirectionEntries} \flowDir{\flowDirection,t})\paramInflowValueUB{v}{t}  \quad\forall v\in\setBoundaryNodes \label{eq:ns_fenceNodes_flowDirection2}
  \end{align}

\paragraph{Flow direction exit pressures}
  Apart from the consequences that the flow direction choice has on the corresponding boundary node inflows, it may also influence the upper pressure bounds on some boundary nodes.
  Each node $v\in\setBoundaryNodesExitP\subseteq\setBoundaryNodes$ is given an upper pressure bound \pressUBExitFG{v}, which is only active if the node is in the outflow set of the currently active flow direction, i.e., serving as exit of the network station.
  The corresponding constraint is the following for each time $t\in\setTimestepsNoZero$
  \begin{align}
    \pressI{v,t} &\leq \pressUBExitFG{v} + ( 1 - \sum_{\flowDirection=(\flowDirectionEntries, \flowDirectionExits)\in\setFlowDirectionsNavi : v\in\flowDirectionExits} \flowDir{\flowDirection,t})(\pressUB{v,t} - \pressUBExitFG{v}) \quad\forall v\in\setBoundaryNodesExitP \label{eq:ns_fenceNodes_exitPressure}
  \end{align}

\paragraph{Flow direction conditions}
  As the final constraints concerning flow directions, there exist in some network stations a special set of conditions \setFlowDirectionConditions, which concern the amount of flow over sets of boundary nodes.
  These conditions must be met for the flow direction which they are associated with to be active.
  Each condition $w = (\flowDirection,\setVertices^{w_1},\setVertices^{w_2})\in\setFlowDirectionConditions$ states that the flow over a set of boundary nodes $\setVertices^{w_1}$ has to be smaller than the flow over a second set of boundary nodes $\setVertices^{w_2}$ if $f$ is selected.
  Note, that the flow over a set of boundary nodes $\setVertices^{w}$ for time $t$ is defined as $\sum_{v\in\setVertices^{w}} |\inflowI{v,t}|$, which is potentially a non-linear expression due to the absolute value.
  However, each boundary node set $\setVertices^{w}$ which is part of some $w\in\setFlowDirectionConditions$ is known to be either a subset of \flowDirectionEntries or a subset of \flowDirectionExits of the corresponding flow direction $\flowDirection=(\flowDirectionEntries, \flowDirectionExits)$.
  For this reason, we always know the sign of the flow over $\setVertices^{w}$ in advance and hence using the following function definition for easier notation
  \begin{equation*}
    \text{sgn}: \setFlowDirectionsNavi \times \setBoundaryNodes \rightarrow \{-1,1\} \,,\, \left((\flowDirectionEntries, \flowDirectionExits),v\right) \rightarrow
    \begin{cases}
      1 & \text{ if } v \in \flowDirectionEntries \\
     -1 & \text{ if } v \in \flowDirectionExits
    \end{cases}
  \end{equation*}
  we can write the constraints for each $t \in \setTimestepsNoZero$ as
  \begin{align}
    \sum_{v \in \setVertices^{w_1}} \text{sgn}(f,v) \inflowI{v,t} &- \sum_{v \in \setVertices^{w_2}} \text{sgn}(f,v) \inflowI{v,t} \notag\\
    &\leq (1 - \flowDir{f,t})C_1 \quad & \forall (f,\setVertices^{w_1},\setVertices^{w_2}) \in \setFlowDirectionConditions. \label{eq:flowDir_condition}
  \end{align}

  Here $C_1$ denotes a big-M constant, which can be set as follows:
  \begin{align*}
    C_1 =&  \sum_{v \in \setVertices^{w_1}:v\in\flowDirectionEntries} \max(0,\paramInflowValueUB{v}{t}) - \sum_{v \in \setVertices^{w_1}:v\in\flowDirectionExits} \min(0,\paramInflowValueLB{v}{t}) \\
    &- \left( \sum_{v \in \setVertices^{w_2}:v\in\flowDirectionEntries} \max(0,\paramInflowValueLB{v}{t}) - \sum_{v \in \setVertices^{w_2}:v\in\flowDirectionExits} \min(0,\paramInflowValueUB{v}{t}) \right)
  \end{align*}

\subsection{Scenario and initial state}
\label{sec:scenario_initState}

For the future we are given scenario values for the boundaries of the station in terms of pressure and inflow.
While we are given one pressure value \demandPressure{v,t} per boundary node $v\in\setBoundaryNodes$ for each future time point $t\in\setTimestepsNoZero$, the flow demands are only given for sets of boundary nodes, which are called the \emph{fence groups} of the network station forming the set \setNSFenceGroups.
For each set $g\in\setNSFenceGroups$, which can also consist of only a single boundary node, and each future time point $t\in\setTimestepsNoZero$, the sum of inflows should be equal to the given demand value \demandInflow{g,t}.

We do not require strict obedience of the given demand values \demandPressure{v,t} and \demandInflow{g,t}, but instead allow deviations from them, which will be punished in the objective function.
These deviations are captured in the \emph{slack variables} \slack.
For pressure we have \slackPPos{v,t} and \slackPNeg{v,t} which capture the positive and respectively negative difference between the pressure values of boundary node $v$ at future time $t$ and the given demand \demandPressure{v,t}.
We then have the slack variables \slackQPos{v,t} and \slackQNeg{v,t} associated with the inflow,  which capture the positive and respectively negative contribution to the difference between the inflow demand \demandInflow{g,t} of fence group $g\in\setNSFenceGroups$ at each future time step $t$ and the sum of the corresponding inflow variables of each boundary node $v$ in the fence group $g$.

The described relations can be modeled for each future time step $t\in\setTimestepsNoZero$ as:
\begin{align}
  \demandPressure{v,t} &= \pressI{v,t} - \slackPPos{v,t} + \slackPNeg{v,t} \quad \forall v\in\setBoundaryNodes \label{eq:slack_pressure}\\
  \demandInflow{g,t} &= \sum_{v\in g}\left( \inflowI{v,t} - \slackQPos{v,t} + \slackQNeg{v,t} \right) \quad \forall g\in\setNSFenceGroups \label{eq:slack_flow}
\end{align}
Note that by using inflow slacks, it is possible to choose a flow direction for the station that does not fit the inflow demand values \demandInflow{g,t} in terms of the constraints~\eqref{eq:ns_fenceNodes_flowDirection1}-\eqref{eq:ns_fenceNodes_flowDirection2}.

The second set of prescribed values are those of the initial state.
A complete list of given values are as follows:
The initial pressures $\pressI{v,0}$ for each node $v\in\setVertices$,
the in- and outflow values \mFlowI{l,a,0} and \mFlowI{r,a,0} for each pipe $(l,r)=a\in\setPipes$,
the flow values \mFlowI{a,0} for each non-pipe arc $a\in\setArcs\setminus\setPipes$,
the operation mode of the network station fixing \opMode{o,0} for each $o\in\setNSModes$,
the thereby determined values of the variables
\modeOp{a,0}, \modeBy{a,0}, \modeCl{a,0}, \modeAc{a,0} for each corresponding valve and compressor station arc $a\in\setValves\cup\setCompressorStations$,
the values of \cfg{c,a,0} for each compressor station $a\in\setCompressorStations$ and corresponding configuration $c\in\setCompressorConfigurations{a}$,
and finally the modes \modeBy{a,0}, \modeCl{a,0}, \modeAc{a,0} of all regulators $a\in\setRegulators$.
All these variables are actually parameters of the model and fixed to the corresponding value.

\subsection{Objective}
\label{sec:objective}

As already describe above, the objective function should punish deviation from the given future scenario, while simultaneously favoring those solutions with a stable control of the single elements.
While the first part can easily be described by using the slack variables introduced in Section~\ref{sec:scenario_initState}, we still need to define a measure for the stability.
To do so, we first quantify the discrete changes of the control in binary variables, i.e., the change of the network station into a new operation mode at time $t$ in variable \opModeChg{t} as well as the change into a new mode of regulator $a$ at time $t$ in variable \regModeChg{a,t}.
Furthermore, we capture the start of compressor unit $u$ at time $t$ in the variable \unitStart{u,t}, since starting a compressor is a very time and energy intensive action and should therefore be avoided if possible.
The mode changes of valves and compressor stations are not tracked separately, as the start-up of compressor stations is already penalized and the change of valves only happens for multiple elements at once in the context of operation mode changes, see the Equation~\eqref{eq:valve_opMode_coupling}.
This described variable behavior is realized by the following constraints for each future time step $t\in\setTimestepsNoZero$, where we denote the set of configurations of the containing compressor station $a\in\setCompressorStations$ which use compressor unit $u\in\setCompressorUnitsI{a}$ by \setCompressorConfigurationsUsingUnit{u}{a}.
\begin{align}
  \opModeChg{t} &\geq \opMode{o,t} - \opMode{o,t-1} \quad \forall o\in\setNSModes \label{eq:opMode_change_definition1}\\
  \opModeChg{t} &\leq 2 - \opMode{o,t} - \opMode{o,t-1} \quad \forall o\in\setNSModes \label{eq:opMode_change_definition2}\\
  \regModeChg{a,t} &\geq \modeCl{a,t} - \modeCl{a,t-1} \quad \forall a\in\setRegulators \label{eq:regulator_mode_change_definition_first}\\
  \regModeChg{a,t} &\leq 2 - \modeCl{a,t} - \modeCl{a,t-1} \quad \forall a\in\setRegulators\\
  \regModeChg{a,t} &\geq \modeBy{a,t} - \modeBy{a,t-1} \quad \forall a\in\setRegulators\\
  \regModeChg{a,t} &\leq 2 - \modeBy{a,t} - \modeBy{a,t-1} \quad \forall a\in\setRegulators\\
  \regModeChg{a,t} &\geq \modeAc{a,t} - \modeAc{a,t-1} \quad \forall a\in\setRegulators\\
  \regModeChg{a,t} &\leq 2 - \modeAc{a,t} - \modeAc{a,t-1} \quad \forall a\in\setRegulators \label{eq:regulator_mode_change_definition_last}\\
  \unitStart{u,t} &\geq \sum_{c\in\setCompressorConfigurationsUsingUnit{u}{a}} \cfg{c,a,t} - \sum_{c\in\setCompressorConfigurationsUsingUnit{u}{a}} \cfg{c,a,t-1} \quad \forall u\in\setCompressorUnitsI{a}\forall a\in\setCompressorStations \label{eq:unit_start_definition}
\end{align}

In order to obtain a smooth operation for network situations without discrete mode switching, we add variables tracking the change of the operation point of single elements, i.e., their corresponding changes in flow, incoming pressure and outgoing pressure.
We do this for all elements with an actual operation point, i.e., regulators and compressor stations in active mode, while ignoring times in which the network station operation mode or regulator mode has just been changed.
The variables \rgPLChg{a,t}, \rgPRChg{a,t}, \rgQChg{a,t} representing changes of the incoming pressure, outgoing pressure and flow of an active regulator $a$ respectively \csPLChg{a,t}, \csPRChg{a,t}, \csQChg{a,t} representing the corresponding value changes of an active compressor station $a$ can be established using the following constraints for each $(l,r)=a\in\setRegulators$ and each $t\in\setTimestepsNoZero$
\begin{align}
  \pressI{l,t} - \pressI{l,t-1} &\leq \rgPLChg{a,t} + (\modeBy{a,t} + \modeCl{a,t} + \regModeChg{a,t})(\pressUB{l,t}-\pressLB{l,t-1}) \label{eq:regulator_operating_point_first}\\
  \pressI{l,t-1} - \pressI{l,t} &\leq \rgPLChg{a,t} + (\modeBy{a,t} + \modeCl{a,t} + \regModeChg{a,t})(\pressUB{l,t-1}-\pressLB{l,t})\\
  \pressI{r,t} - \pressI{r,t-1} &\leq \rgPRChg{a,t} + (\modeBy{a,t} + \modeCl{a,t} + \regModeChg{a,t})(\pressUB{r,t}-\pressLB{r,t-1})\\
  \pressI{r,t-1} - \pressI{r,t} &\leq \rgPRChg{a,t} + (\modeBy{a,t} + \modeCl{a,t} + \regModeChg{a,t})(\pressUB{r,t-1}-\pressLB{r,t})\\
  \mFlowI{a,t} - \mFlowI{a,t-1} &\leq \rgQChg{a,t} + (\modeBy{a,t} + \modeCl{a,t} + \regModeChg{a,t})(\mFlowUB{a,t}-\mFlowLB{a,t-1})\\
  \mFlowI{a,t-1} - \mFlowI{a,t} &\leq \rgQChg{a,t} + (\modeBy{a,t} + \modeCl{a,t} + \regModeChg{a,t})(\mFlowUB{a,t-1}-\mFlowLB{a,t}) \label{eq:regulator_operating_point_last}
\end{align}
respectively for each $(l,r)=a\in\setCompressorStations$ and each $t\in\setTimestepsNoZero$
\begin{align}
  \pressI{l,t} - \pressI{l,t-1} &\leq \csPLChg{a,t} + (\modeBy{a,t} + \modeCl{a,t} + \opModeChg{t})(\pressUB{l,t}-\pressLB{l,t-1}) \label{eq:compressorStation_operating_point_first}\\
  \pressI{l,t-1} - \pressI{l,t} &\leq \csPLChg{a,t} + (\modeBy{a,t} + \modeCl{a,t} + \opModeChg{t})(\pressUB{l,t-1}-\pressLB{l,t})\\
  \pressI{r,t} - \pressI{r,t-1} &\leq \csPRChg{a,t} + (\modeBy{a,t} + \modeCl{a,t} + \opModeChg{t})(\pressUB{r,t}-\pressLB{r,t-1})\\
  \pressI{r,t-1} - \pressI{r,t} &\leq \csPRChg{a,t} + (\modeBy{a,t} + \modeCl{a,t} + \opModeChg{t})(\pressUB{r,t-1}-\pressLB{r,t})\\
  \mFlowI{a,t} - \mFlowI{a,t-1} &\leq \csQChg{a,t} + (\modeBy{a,t} + \modeCl{a,t} + \opModeChg{t})(\mFlowUB{a,t}-\mFlowLB{a,t-1})\\
  \mFlowI{a,t-1} - \mFlowI{a,t} &\leq \csQChg{a,t} + (\modeBy{a,t} + \modeCl{a,t} + \opModeChg{t})(\mFlowUB{a,t-1}-\mFlowLB{a,t}). \label{eq:compressorStation_operating_point_last}
\end{align}
Note that we needed to define the upper bound constraints for the discrete change variables \opModeChg{t} and \regModeChg{a,t} to allow them to be equal to one only if there really is a discrete change.
Otherwise it might have been possible to set the change variable to 1 although there is no actual discrete change and thereby avoid high costs imposed by the continuous change variables, which is not a desired behavior.

Finally, we are able to state our objective function, which minimizes the weighted sum of the change variables and the slack variables defined in Section~\ref{sec:scenario_initState} as
\begin{align}
  \min &\text{ objective}:=\notag\\
  \sum_{t\in\setTimestepsNoZero} \bigg(
   & \big(\granularityI{t}-\granularityI{t-1}\big)
   \sum_{v\in\setBoundaryNodes} \big(\costSlackP \cdot (\slackPPos{v,t} + \slackPNeg{v,t}) + \costSlackQ \cdot (\slackQPos{v,t} + \slackQNeg{v,t})\big) \notag\\
  &+ \costOpModes\cdot\opModeChg{t} \notag\\
  &+ \sum_{a\in\setRegulators} \costRgModeChg\cdot\regModeChg{a,t} \label{eq:objective}\\
  &+ \sum_{u\in\setCompressorUnitsI{a}, a\in\setCompressorStations} \costUnitStarts\cdot\unitStart{u,t} \notag\\
  &+ \sum_{a\in\setRegulators} \big( \costRgPLChg\cdot\rgPLChg{a,t} + \costRgPRChg\cdot\rgPRChg{a,t} + \costRgQChg\cdot\rgQChg{a,t} \big)\notag\\
  &+ \sum_{a\in\setCompressorStations} \big( \costCsPLChg\cdot\csPLChg{a,t} + \costCsPRChg\cdot\csPRChg{a,t} + \costCsQChg\cdot\csQChg{a,t} \big)\bigg), \notag
\end{align}
where the $\cost^{*}$ parameters denote the corresponding positive weights given to the single quantities.
Note that the weights for pressure and flow slack variables are additionally multiplied by the length of the corresponding time interval.
We do this, since the slacks track deviation amounts over time and we therefore need to take the time interval length into account for their objective function value coefficients.

\subsection{Final model}
\label{sec:final_model}

Putting everything together, we can formulate our problem in the following transient gas flow model \MIP:
\begin{align*}
  \min\quad & \eqref{eq:objective} \\
  \text{s.t.} \quad \forall t\in\setTimestepsNoZero \qquad
  & \eqref{eq:compressorStation_OneModeOrConfig}-\eqref{eq:compressorStation_last1}, \eqref{eq:compressorStation_first2}-\eqref{eq:compressorStation_last2}, \eqref{eq:cs_bypass_opMode_coupling}-\eqref{eq:cs_closed_opMode_coupling}, \eqref{eq:compressorStation_operating_point_first}-\eqref{eq:compressorStation_operating_point_last} \quad \forall a\in\setCompressorStations \\
  & \eqref{eq:compressorStation_cfg_first}-\eqref{eq:compressorStation_cfg_last}, \eqref{eq:cs_cfg_opMode_coupling} \quad \forall a\in\setCompressorStations \quad \forall c\in\setCompressorConfigurations{a} \\
  & \eqref{eq:compressorStation_cfg_facets} \quad \forall a\in\setCompressorStations \quad \forall c\in\setCompressorConfigurations{a} \quad \forall (w,x,y,z)\in\setConfigurationFacets{c} \\
  & \eqref{eq:pipes_constVelo_continuity}-\eqref{eq:pipes_constVelo_momentum} \quad \forall a\in\setPipes \\
  & \eqref{eq:resistors_constVelo}\quad \forall a\in\setResistors \\
  & \eqref{eq:valves_first}-\eqref{eq:valves_last},\eqref{eq:valve_opMode_coupling} \quad \forall a\in\setValves \\
  \MIP \rule{2.5cm}{0.0mm}
  & \eqref{eq:regulator_mode}-\eqref{eq:regulator_mode_definition_last},\eqref{eq:regulator_mode_change_definition_first}-\eqref{eq:regulator_mode_change_definition_last},\eqref{eq:regulator_operating_point_first}-\eqref{eq:regulator_operating_point_last} \quad \forall a \in \setRegulators  \\
  & \eqref{eq:flowConservation_boundaryNodes},\eqref{eq:ns_fenceNodes_flowDirection1}-\eqref{eq:ns_fenceNodes_flowDirection2},\eqref{eq:slack_pressure} \quad \forall v\in\setBoundaryNodes \\
  & \eqref{eq:flowConservation_innerNodes} \quad \forall v\in\setInnerNodes \\
  & \eqref{eq:opMode_choice}, \eqref{eq:flowDir_choice} \\
  & \eqref{eq:ns_om_fd_coupling},\eqref{eq:opMode_change_definition1}-\eqref{eq:opMode_change_definition2} \quad \forall o\in\setNSModes \\
  & \eqref{eq:ns_fenceNodes_exitPressure} \quad \forall v\in\setBoundaryNodesExitP \\
  & \eqref{eq:flowDir_condition} \quad \forall w\in\setFlowDirectionConditions \\
  & \eqref{eq:slack_flow} \quad \forall g\in\setNSFenceGroups \\
  & \eqref{eq:unit_start_definition} \quad \forall a\in\setCompressorStations \quad u\in\setCompressorUnitsI{a}
\end{align*}
Note that we apply the constraints only starting from time step 1 explicitly excluding the initial time step 0.
We do this since the initial pressure and flow values as well as the initial modes described in Section~\ref{sec:scenario_initState} are not guaranteed to fit our model, and simply serve as a starting point for the calculations.

\section{Feasible operating range of compressor station configurations}
\label{sec:compressor_configurations_definition}

As already explained in Section~\ref{sec:compressor_station} each compressor station arc $(l,r)=a\in\setCompressorStations$ has an inherent substructure.
It represents a set of compressor units \setCompressorUnitsI{a}, which are the actual compressing elements and can be combined in a serial and/or parallel fashion, to either allow for a higher compression ratio, a higher flow rate or a mixture of both.
The set of all feasible serial-parallel combinations is called the set of configurations \setCompressorConfigurations{a} of a compressor station.
For each of these configurations, we will describe in this section how to obtain their polytope description given as intersection of the set of half spaces \setConfigurationFacets{c}, which is used in the model as described in Section~\ref{sec:compressor_station}.
The polytope of a configuration is created based on the polytopes of the compressor units, which are created from combining the corresponding feasible operation ranges with a maximum power restriction.
Note that we drop the time index in this section for the ease of notation.

\subsection{Feasible operating range for a single compressor unit}

A compressor unit is a combination of a single compressor machine (or just a compressor), which increases the gas pressure in flow direction, and a corresponding drive, providing the power needed to run the compressor.
For each compressor machine we are given a feasible operating range as polytope in the space $(\frac{\pressI{r}}{\pressI{l}}, \vFlow)$, where the volumetric flow rate \vFlow is given as \begin{equation*}
  \vFlow = \mFlow / \dens_l\text {, \quad with } \dens_l = \frac{\pressI{l}}{\sGasConst \temp \zFactI{l}}.
\end{equation*}
Note, that since we explicitly consider density and pressure at the incoming node $l$ of the compressor, we use the compressibility factor $\zFactI{l}:=\zFact(\pressI{l,0})$ instead of $\zFactI{a}$, here as well as in the rest of this Section.

Usually, the feasible operating range, sometimes also called ``characteristic diagram'' or ``performance curve'', is given as area in the dimensions $(\adEnth, \vFlow)$ restricted by a set of possibly concave quadratic curves, see e.g.~\cite{FugGeiGolMor2015}\cite{OdoMus2009}.
The quantity \adEnth denotes the \emph{specific change in adiabatic enthalpy} and is defined as
\begin{equation}
  \adEnth = \sGasConst \temp \zFactI{l} \frac{\isenExp}{\isenExp-1} \left[ \left( \frac{\pressI{r}}{\pressI{l}} \right)^{\frac{\isenExp-1}{\isenExp}} - 1 \right], \label{eq:Had_definition}
\end{equation}
using for the \emph{isentropic exponent} \isenExp the constant value $1.296$, as stated in~\cite{FugGeiGolMor2015}.
The transformation of such a feasible operating range using \adEnth into our format is easily doable, since there is a unique transformation from \adEnth to $\frac{\pressI{r}}{\pressI{l}}$ obtained by simply rearranging \eqref{eq:Had_definition}.
The diagram then just has to be linearized by approximation or relaxation to obtain the polytope description in the desired space.

\begin{figure}[ht]
  \centering

  \begin{subfigure}[c]{0.48\textwidth}
  \includegraphics[trim={1.0cm 1.8cm 1.0cm 1.0cm},clip,width=\textwidth]{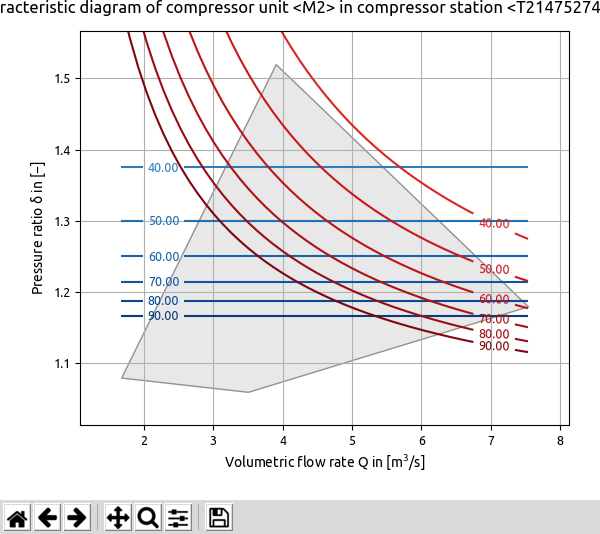}
  \subcaption{Original power bound}
  \label{fig:unit_charDiag}
  \end{subfigure}
  \begin{subfigure}[c]{0.48\textwidth}
  \includegraphics[trim={1.0cm 1.8cm 1.0cm 1.0cm},clip,width=\textwidth]{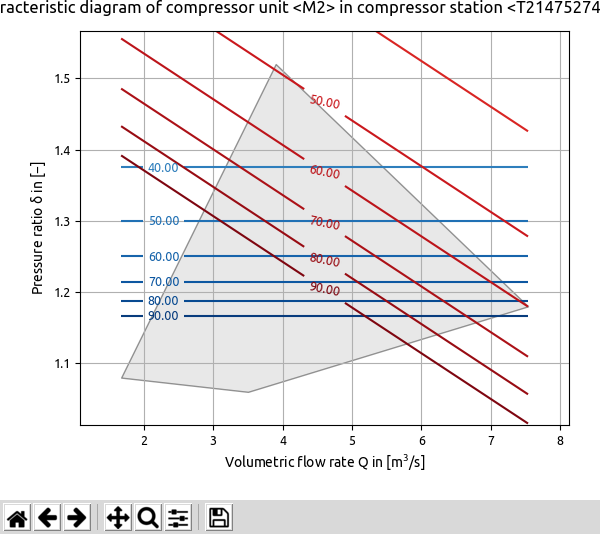}
  \subcaption{Linearized power bound}
  \label{fig:unit_charDiag_linearized}
  \end{subfigure}
  \caption{The feasible operating range of a compressor unit.
  The grey region shows the operating range given as a polytope.
  For different values of the incoming pressure \pressI{l}
  the blue lines represent the upper bound on the absolute pressure increase~\machinePDiffUB and the red lines illustrate the power bound $\ub{\power}$.
  While the left picture shows the original non-linear non-convex power bound, the right pictures shows the linearized version, see Section~\ref{sec:powerBound_linearization}}.
  \label{fig:unit_charDiag_both}
\end{figure}

In addition to the feasible operating range polytope, each compressor machine is given an upper bound on the absolute pressure increase $\machinePDiffUB \geq \pressI{r}-\pressI{l}$ and an upper bound on the maximum power to use $\ub{\power}$ based on the power of the compressor drive.
The power needed for compression depends on the above defined \adEnth as well as the mass flow and is given as
\begin{equation}
  \power = \frac{\mFlow \adEnth}{\adEff}
  = \frac{\mFlow}{\adEff}\sGasConst \temp \zFactI{l} \frac{\isenExp}{\isenExp-1} \left[ \left( \frac{\pressI{r}}{\pressI{l}} \right)^{\frac{\isenExp-1}{\isenExp}} - 1 \right].
  \label{eq:power_definition}
\end{equation}
Here \adEff denotes the \emph{adiabatic efficiency} of the compression, which in theory depends on the actual point of operation in the feasible region but is here assumed to be a given constant per compressor unit.
An example of a feasible operating range of a compressor unit is given in Figure~\ref{fig:unit_charDiag}, where different levels of the maximum pressure difference bound and the maximum power bound are given based on different values for the incoming pressure $\pressI{l}$.

\begin{figure}[ht]
  \centering

  \begin{subfigure}[c]{0.48\textwidth}
  \includegraphics[trim={1.0cm 1.8cm 1.0cm 1.0cm},clip,width=\textwidth]{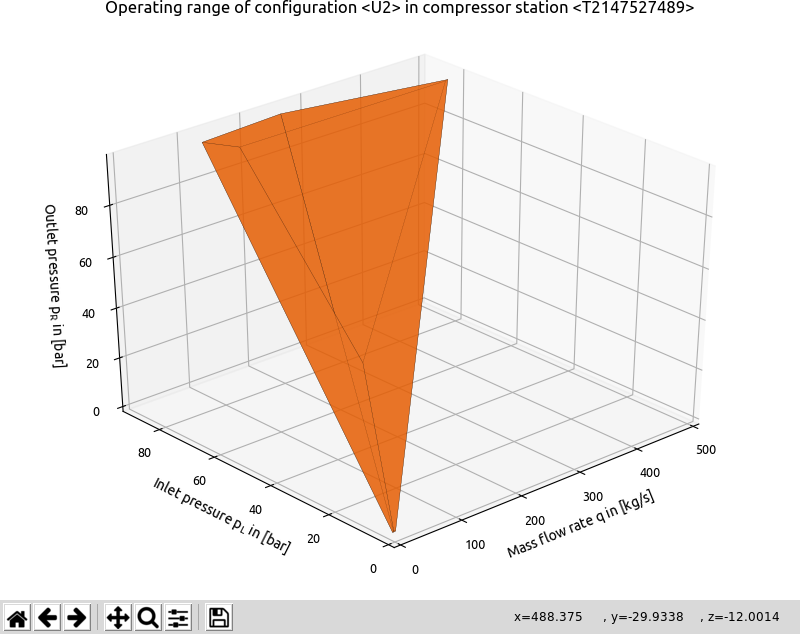}
  \subcaption{Without power bound}
  \label{fig:unit_3D_noPowerBound}
  \end{subfigure}
  \begin{subfigure}[c]{0.48\textwidth}
  \includegraphics[trim={1.0cm 1.8cm 1.0cm 1.0cm},clip,width=\textwidth]{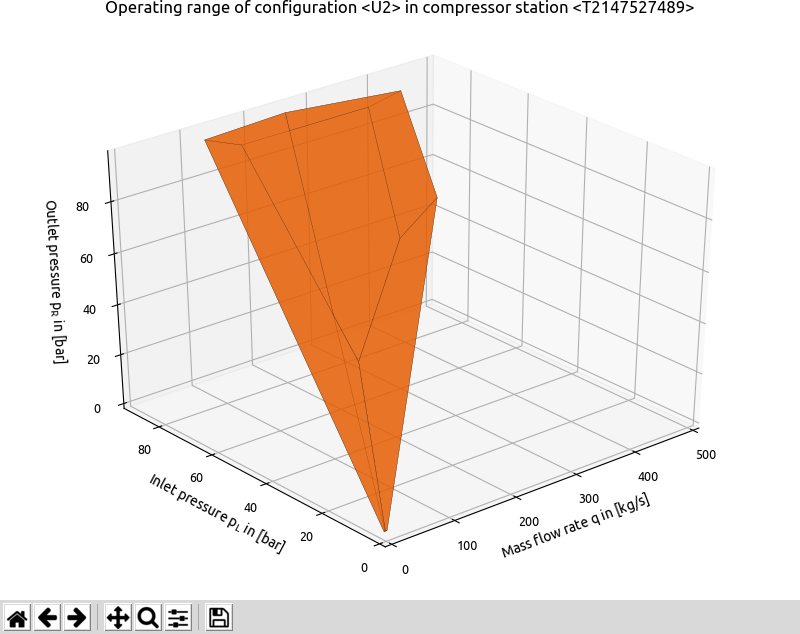}
  \subcaption{With power bound}
  \label{fig:unit_3D_PowerBound}
  \end{subfigure}
  \caption{The feasible operating range of a compressor unit in the space $(\pressI{l}, \pressI{r}, \mFlow)$, computed from the two dimensional operating range shown in Figure~\ref{fig:unit_charDiag_both}.
  While the left picture shows the lifted polytope based on the original feasible operating range, the maximum absolute pressure difference and the end nodes pressure bounds, the right pictures also includes the linearized power bound, see Section~\ref{sec:powerBound_linearization}}.
\end{figure}

Ignoring the power bound for a moment, we will now lift the feasible operating range into the $(\pressI{l}, \pressI{r}, \mFlow)$ space we are interested in.
Therefore, we first transform each of the faces $a_0 + a_1 \vFlow + a_2 \frac{\pressI{r}}{\pressI{l}} \leq 0$ of the original polytope into constraints  $\tilde{a}_0 \pressI{l} + \tilde{a}_1 \pressI{r} + \tilde{a}_2 \mFlow \leq 0$ of the higher dimensional space using the equation of state for real gases \eqref{eq:realGasState}:
\begin{align*}
  a_0 + a_1 \vFlow + a_2 \frac{\pressI{r}}{\pressI{l}} &\leq 0 \\
  \quad \Leftrightarrow \quad
  a_0 + a_1 \frac{\mFlow}{\rho_L} + a_2 \frac{\pressI{r}}{\pressI{l}} &\leq 0 \\
  \quad \Leftrightarrow \quad
  a_0 + a_1 \frac{\mFlow \sGasConst \temp \zFactI{l}}{\pressI{l}} + a_2 \frac{\pressI{r}}{\pressI{l}} &\leq 0 \\
  \quad \Leftrightarrow \quad
  a_0 \pressI{l} + a_1 \sGasConst \temp \zFactI{l} \mFlow + a_2 \pressI{r} &\leq 0 \\
  \quad \Leftrightarrow \quad
  \tilde{a}_0 \pressI{l} + \tilde{a}_1 \pressI{r} + \tilde{a}_2 \mFlow &\leq 0.
\end{align*}
To bound the polyhedron described by the new constraints, we add the restriction of the absolute pressure difference as well as two pressure bounds of the end nodes of the compressor station arc $(l,r)=a\in\setCompressorStations$ containing this machine
\begin{align*}
  \pressI{r} -\pressI{l} &\leq \machinePDiffUB \\
  \pressI{l} &\geq \lb{p}_l \\
  \pressI{r} &\leq \ub{p}_r.
\end{align*}
A picture of the three dimensional polytope resulting from the feasible operating range of Figure~\ref{fig:unit_charDiag_both} can be seen in Figure~\ref{fig:unit_3D_noPowerBound}.

\subsection{Power bound linearization}
\label{sec:powerBound_linearization}

Until now, we have ignored the maximum power constraint in the three dimensional feasible operation range polytope, that restricts the value of the available power \power for compression.
Figure~\ref{fig:unit_charDiag} shows the constraint~$ \power \leq \ub{\power} $ cuts into the original two dimensional feasible operating range in a non-linear and non-convex fashion.
The same holds for the feasible operating range representation in $(\pressI{l},\pressI{r},\mFlow)$.
In the following, we are going to derive a linear approximation to this constraint that can then be
added to the operating range polytope.

Therefore, we generate a set of $N$ random sample points from within the three dimensional operating range polytope, which we represent by the vectors $ \mathbf{\pressI{l}}, \mathbf{\pressI{r}}, \mathbf{\mFlow} \in \mathbb{R}^N$.
We used a rejection-free sampling method, which is based on the sampling from a 3D tetrahedron, see~\cite{RoccCig2001}.
The idea behind this method is to sample uniformly from a through the tetrahedron generated parallelepiped and then project the sample back to the tetrahedron.
To extend this method from tetrahedra to general polytopes, we first compute a triangulation of the polytope.
We compute the volumes of the tetrahedra and use their relative volume as a probability function.
This allows us to first randomly select a tetrahedron and then sample from it using the in~\cite{RoccCig2001} described method.
Our approach has high preprocessing costs of computing the triangulation tetrahedra and their volume but then it is very cheap to generate samples.
On the contrary for rejection based sampling there is no lower bound for how many samples we need to generate for a single sample point, since the polytope can be arbitrarily small in comparison to the space we sample from (typically the enclosing axis-parallel cube).
For the operating range polytopes we observed a significant speed-up while generating 50000 samples.

For each point, we then determine the corresponding compressor power using Equation~\eqref{eq:power_definition} and store these values again in a vector $ \mathbf{\power} \in \mathbb{R}^N$.
The goal is now to obtain a linear approximation of the power function, i.e., an approximation of the form
\begin{align*}
  \mathbf{\power} \approx a_0 + a_1 \mathbf{\pressI{l}} + a_2 \mathbf{\pressI{r}} + a_3 \mathbf{\mFlow}.
\end{align*}
We achieve this by applying an ordinary least-squares method in order to determine the coefficients of the linear function as
\begin{align*}
  \min\limits_{a_0,a_1,a_2,a_3} \mid\mid \mathbf{\power} - ( a_0 + a_1 \mathbf{\pressI{l}} + a_2 \mathbf{\pressI{r}} + a_3 \mathbf{\mFlow} ) \mid\mid^2.
\end{align*}
Finally, we can formulate the new linearized power bound constraint based on the obtained solution
values $ (\ub{a}_0,\ub{a}_1,\ub{a}_2,\ub{a}_3) $ as
\begin{align*}
  \ub{a}_0 + \ub{a}_1 \pressI{l} + \ub{a}_2 \pressI{r} + \ub{a}_3 \mFlow \leq \ub{\power},
\end{align*}
which is then added to the three dimensional polytope to create the final three description of the feasible operating range of a compressor unit.
An example the final polytope is illustrated in Figure~\ref{fig:unit_3D_PowerBound}, while the linearized power bound projected to the original two dimensional operating range can be seen in Figure~\ref{fig:unit_charDiag_linearized}.

\subsection{Feasible operating range for a compressor station configuration}

As the final step, we will create the polytope description for each configuration by combining the polytopes of the used compressor units.
The procedure was originally described in~\cite{HumFugHilKoc2015}, while we exactly follow the steps of the variant described in~\cite{WalHil2017}.

Each configuration $c$ is given as a serial sequence $s_1, \dots, s_{n_c}$ of parallel compressor machine arrangements, where combining compressors in series allows for higher output pressures by multi-step compression while parallel compression increases the throughput in terms of flow.
We call such a parallel machine arrangement \emph{stage} and denote by \setCompressorUnitsI{s} the set of compressor units combined in stage $s\in\{s_1, \dots, s_{n_c}\}$.

We will now start with the definition of the feasible operation range polytope $P_s$ of such a stage $s\in\{s_1, \dots, s_{n_c}\}$.
For each compressor unit $u\in\setCompressorUnitsI{s}$ the corresponding polytope is denoted by $P_u$, and using this we can describe $P_s$ as
\begin{align*}
  P_s := \{\quad(\pressI{l},\pressI{r},\mFlow) \quad | & \quad \forall u\in\setCompressorUnitsI{s}\quad\exists (\pressI{l,u}, \pressI{r,u}, \mFlowI{u}) \in P_u\\
  \text{ with }
  \pressI{l} &= \pressI{l,u'} \quad \forall u'\in\setCompressorUnitsI{s}, \\
  \pressI{r} &= \pressI{r,u'} \quad \forall u'\in\setCompressorUnitsI{s}, \\
  \mFlow &= \sum_{u'\in\setCompressorUnitsI{s}} \mFlowI{u'} \quad \}.
\end{align*}
In other words, a valid operation point of the stage is represented by each unit operating at the same incoming and outgoing pressures, while the the mass flow through the stage is the sum of mass flow through the individual units.

In a similar fashion, we can then define the polytope $P_c$ for the overall configuration.
Here, the mass flow through all the stages stays the same, while the outgoing pressure of some stage $s_i$ in the sequence has to match the incoming pressure of the subsequent stage $s_{i+1}$.
Using this logic $P_c$ can be defined as
\begin{align*}
  P_c := \{\quad(\pressI{l},\pressI{r},\mFlow) \quad | & \quad \forall s\in\{s_1, \dots, s_{n_c}\}\quad\exists (\pressI{l,s}, \pressI{r,s}, \mFlowI{s}) \in P_s\\
  \text{ with }
  \pressI{l} &= \pressI{l,s_1}, \\
  \pressI{r} &= \pressI{r,s_{n_c}}, \\
  \pressI{r,s_i} &= \pressI{l,s_{i+1}} \quad \forall i\in\{1,\dots,n_c - 1\}, \\
  \mFlow &= \mFlowI{s'} \quad \forall s'\in\{s_1, \dots, s_{n_c}\} \quad \}.
\end{align*}
Finally, the set of half spaces \setConfigurationFacets{c} used in Section~\ref{sec:compressor_station} to define the feasible operating range of the configuration $c$ are simply the facets of the polytope $P_c$ representing the feasible operating range.

Note, that due to the symmetric polytope creation the parallel compressor units of a stage do not need a specific order.
In contrast, the serial stage sequence is indeed important since in general two sequences using the same stages but in a different order will have different operating ranges.

\section{Specialized network station algorithm}
\label{sec:stationModelAlgorithm}

  The MIP model \MIP presented in Section~\ref{sec:final_model} turns out to be quite challenging according to our experiments even though we are only considering parts of a larger gas network by restricting ourselves to network stations.
  We therefore created a specialized algorithm to solve the problem for network stations.
  The baseline insight of it is that the elements in network stations are all very close to each other, making the corresponding pipes inside the station relatively short, see also Table~\ref{tab:stationStatistics} in the computational Section~\ref{sec:computational_results}.
  This means that their capability to store gas, which is often referred to as \emph{linepack}, is insignificant in comparison to the long pipelines in between the network stations.
  Thus there is no possibility to ``prepare for the future'', i.e., use linepack to pre-transport gas to handle upcoming critical demand situations inside the network station itself.

  This led us to the idea of splitting the time coupled model \MIP into individual stationary models to determine the best operation mode and flow direction for each individual time step.
  We then use this information to find a solution in terms of operation modes and flow directions over the whole time horizon.
  This solution will respect the transition time conditions, which have not been modeled in \MIP explicitly, see the corresponding part of Section~\ref{sec:gasNetworkStationDefinition}.

  Since our goal is to find a feasible solution for the presented model \MIP of Section~\ref{sec:final_model}, we will not only have to determine the network station operation modes and flow directions, but also all other involved quantities.
  Hence, after determining these values we still need to calculate a transient version of \MIP, which for example also takes care of minimizing the differences in the operation points of the single elements, see Section~\ref{sec:objective}.
  We do so by prescribing the operation modes and flow directions determined in the stationary calculations to make the model tractable.
  The whole algorithm is summarized in Algorithm~\ref{algo:completeStationModel}, in which the creation of the operation mode sequence is already split into the two steps described in Section~\ref{sec:determine_operation_modes}.

\begin{algorithm}[ht]
  \DontPrintSemicolon
  $\var{opModes}, \var{flowDirs} \gets \func{initialSolutionCreation}()$\tcp*[f]{see Section~\ref{sec:initialSolution}}\;
  $\var{opModes}, \var{flowDirs} \gets \func{improvementHeuristic}(\var{opModes}, \var{flowDirs})$\tcp*[f]{see Section~\ref{sec:localImprovement}}\;
  return \func{transientSmoothing}(\var{opModes}, \var{flowDirs})\tcp*[f]{see Section~\ref{sec:smoothing}}\;
  \caption{Specialized network station algorithm}
  \label{algo:completeStationModel}
\end{algorithm}

  Note, that our algorithm does not have the guarantee of finding globally optimal solutions, which would have been the case when directly solving \MIP enhanced by a proper description of the operation mode transition time constraints.
  However, it has the potential to obtain good overall solutions for \MIP, since we assume the operation mode depending objective function weights \costOpModes for changing to a new operation mode and \costUnitStarts for starting a new compressor unit to dominate the other objective weights, see also the used objective function values in our computational experiments stated in Section~\ref{sec:instanceDescription}.

  In this section we will first introduce in \ref{sec:modelVariants} the three different variants of the model \MIP that we will use in our algorithm.
  Afterwards, we will present the algorithm itself, that determines operation modes in \ref{sec:determine_operation_modes} and finish with the so-called smoothing procedure in \ref{sec:smoothing}, in which we will find the feasible transient solution to our problem.

\subsection{Model variants}
\label{sec:modelVariants}

  As mentioned above, we will not solve the complete transient model \MIP directly, but use the following variants of it in our overall solution approach.

\paragraph{\MIPs\, -- Stationary model}

  For the first variant we solve a \underline{s}tationary version of the model to determine the best operation mode and flow direction for one independent time step $t$.
  The necessary changes in \MIP mainly affect the pipe model.
  In the stationary case a pipe no longer has the possibility to store gas, since the incoming and outgoing pipe flows are balanced, i.e., $\mFlowI{l,a,t} = \mFlowI{r,a,t} = \mFlowI{a,t}$ for all pipes $(l,r)=a\in\setPipes$.
  This is due to the Continuity Equation, in which $\partial_t p = 0$ holds as for all time dependent derivatives resulting in the mass flow balance.
  Hence, the Continuity Equation is no longer part of the model and the stationary Momentum Equation~\eqref{eq:pipes_constVelo_momentum} for pipe $(l,r)=a\in\setPipes$ and the stationary time step $t$ we are considering can be stated as
  \begin{equation*}
      \pressI{r,t} - \pressI{l,t}
    + \frac{\fricI{a}\lenI{a}}{4\diamI{a}\areaI{a}}
      \left(\absVelo{l,a} + \absVelo{r,a}\right)\mFlowI{a,t}
    + \frac{\gravAcc\slopeI{a}\lenI{a}}{2\sGasConst\temp\zFactI{a}}
      \left( \pressI{l,t} + \pressI{r,t}\right) = 0.
  \end{equation*}
  Note that we still calculate the fixed velocity based on the initial pressure and flow values from time step 0, since we are looking for a feasible solution for the original model \MIP in the end.

  For all other elements the constraints only consider exactly one time step and we therefore simply apply them for the one time step $t$ of the stationary case.
  The only other part to adjust is the objective function, where we keep the penalties \slackPPos{v,t}, \slackPNeg{v,t}, \slackQPos{v,t}, and \slackQNeg{v,t} for each boundary node $v\in\setBoundaryNodes$ since they are defined based on each individual time step $t$.
  In addition, we will keep tracking the change to a new network station operation mode in variable $\opModeChg{t}$ as well as the start of new compressor units using $\unitStart{u,t}$ for all $u\in\setCompressorUnitsI{a}\,a\in\setCompressorStations$. This is done by calling the model with the parameter \texttt{prevMode}, where \texttt{prevMode} represents the operation mode of the previous time step.
  The other variables tracking changes in regulator modes $\regModeChg{a,t}$ as well as the current point of operation of regulators and compressor stations as $\rgPLChg{a,t}$, $\rgPRChg{a,t}$, $\rgQChg{a,t}$, $\csPLChg{a,t}$, $\csPRChg{a,t}$, and $\csQChg{a,t}$ for all $a\in\setRegulators$ respectively $a\in\setCompressorStations$ will be removed from the model, as well as the constraints~\eqref{eq:regulator_mode_change_definition_first}-\eqref{eq:regulator_mode_change_definition_last} respectively \eqref{eq:regulator_operating_point_first}-\eqref{eq:compressorStation_operating_point_last} defining their behavior.
  The final stationary objective function for the time step $t$ under consideration reads as
\begin{align*}
  \min\quad &  \big(\granularityI{t}-\granularityI{t-1}\big)\cdot\sum_{v\in\setBoundaryNodes} \costSlackP\cdot(\slackPPos{v,t} + \slackPNeg{v,t}) + \costSlackQ\cdot(\slackQPos{v,t} + \slackQNeg{v,t})\\
  +& \, \costOpModes\cdot\opModeChg{t} \\
  +& \sum_{u\in\setCompressorUnitsI{a}, a\in\setCompressorStations} \costUnitStarts\cdot\unitStart{u,t}
\end{align*}

\paragraph{\MIPf\, -- Transient with fixed operation modes and flow directions}

  For this variant, we are given a \underline{f}ixed operation mode $o_t$ and flow direction $\flowDirection_t$ for all future time steps $t\in\setTimestepsNoZero$.
  This is used in our algorithm to finally determine all transient quantities after the decision for an operation mode and a flow direction for each time step has been made.
  By fixing the corresponding variables \opMode{o_t,t} and \flowDir{\flowDirection_t,t} to $1$ for all future time steps $t\in\setTimestepsNoZero$, the majority of binary variables can be replaced by constants, since the operation modes already decide valve and compressor station modes, as well as the configuration of all active compressor stations.
  Only the binary variables \modeAc{a,t}, \modeBy{a,t} and \modeCl{a,t} for the mode of a regulator $a\in\setRegulators$ are still to be decided.

  In addition, a lot of implicating big-M constraints can already be resolved, i.e., we can remove the current formulation and just add the implied constraints, if the corresponding condition is fulfilled.
  Examples for this are the constraints~\eqref{eq:valves_first}-\eqref{eq:valves_last} describing the valve behavior or the flow direction conditions~\eqref{eq:flowDir_condition}.
  Furthermore, we no longer need to use a disjunctive model, but can apply the corresponding constraints~\eqref{eq:compressorStation_define_pL}-\eqref{eq:compressorStation_cfg_facets} of the active mode and/or configuration of time step $t$ directly to the variables \pressI{l,t}, \pressI{r,t} and \mFlowI{a,t} for each compressor station $(l,r)=a\in\setCompressorStations$.

\paragraph{\MIPsf\, -- Stationary with fixed operation mode}

  As a last variant we basically combine the two variants above and use the \underline{s}tationary version of the model with already \underline{f}ixed operation mode.
  Note, that in contrast to model \MIPf the flow direction of the network station is not already given, which results in more binary variables and still to decide big-M constraints.
  However, this variant still results in a very small and rather simple model and we can therefore solve it very often to test the appropriateness of a given operation mode for a certain time step.

\subsection{Determining operation modes}
\label{sec:determine_operation_modes}

  Our algorithm to determine the operation modes of the network station is split into two steps:
  First, we create an initial solution by a greedy, forward oriented procedure presented in \ref{sec:initialSolution}.
  We then in a second step improve this solution by testing if certain operation modes can be replaced by similar ones to find a better sequence of operation modes over time.
  This second step is described in \ref{sec:localImprovement}.

  Regarding the flow directions we will return for each time step $t$ that flow direction, which has been chosen in the optimal solution of the stationary model which determined the returned operation mode for $t$.
  For simplicity of notation, we will not further mention them in the rest of this section.

\subsubsection{Initial solution creation}
\label{sec:initialSolution}

  To find a first feasible sequence of operation modes over time, we follow a rather simple idea.
  In order to keep the number of needed operation mode changes small, we determine an operation mode for time step $t$ by first testing the used operation mode of the previous time step $t-1$ using \MIPsf.
  Only when this previously used operation mode yields a costly solution in terms of the objective function value will we use the general stationary model \MIPs to determine the best operation mode for $t$.
  By this mechanic, we also reduce the amount of calls to the \MIPs model, which are in general much more expensive in terms of computing time then calls to the \MIPsf model.
  A detailed description is given as Algorithm~\ref{algo:initialSolve}.

\begin{algorithm}[ht]
  \DontPrintSemicolon
  \KwData{Operation mode o$_0$ of the network station in the initial state}
  \KwResult{A list of operation modes for each time $\var{t}\in\setTimesteps$}
  $\var{operationModes} \gets \func{list}()$\;
  \var{operationModes}.\func{add}($\var{o}_0$)\;
  \For{$t\in\setTimestepsNoZero$}{
    $\var{oldMode} \gets \var{operationModes}.\func{last}()$\;
    \tcp{call \MIPsf with fixed operation mode oldMode for time t}
    $\var{oldModeFeasible}, \var{oldModeCost}$\nonl\\\quad$\gets \MIPsf(\var{oldMode}, \var{t}, prevMode = \var{oldMode})$\label{line:init_callMIPsf}\;
    \eIf{\var{oldModeFeasible}
    \nonl\algoAnd $\func{modeAvailable}(\var{oldMode},\var{t})$
    \nonl\algoAnd $\var{oldModeCost} < \costOpModes$\label{line:init_checkLastOpMode}}{
      \tcp{operation mode of t-1 is also good for t}
      \var{operationModes}.\func{add}(\var{oldMode})\;
    }
    {
      \tcp{operation mode of t-1 is NOT good for t}
      \tcp{$\Rightarrow$ search best possible valid stationMode for t}
      $\var{validModes} \gets \func{list}()$\label{line:init_lastOpModeNotWorking}\;
      \For{$\var{o}\in\setNSModes$}{
        \If{$\func{modeAvailable}(\var{o},\var{t})$
        \nonl\algoAnd $\func{transitionsWork}(\func{concat}(\var{operationModes}, \func{list}(\var{o})))$
        \nonl\algoAnd $\func{notSoonInfeasible}(\var{t}, newMode = \var{o}, oldMode = \var{oldMode})$}{
          $\var{validModes}.\func{add}(\var{o})$\;
        }
      }
      \tcp{find best from validModes by calling \MIPs for time t}
      $\var{bestMode}, \var{bestModeFeasible}, \var{bestModeCost}$ \nonl\\\quad$\gets \MIPs(\var{t}, \setNSModes = \var{validModes}, prevMode = \var{oldMode})$\label{line:init_callMIPs}\;
      \lIf{not \var{bestModeFeasible}}{abort without solution\label{line:abortInitialSolve}}
      \tcp{bestMode is best choice for t}
      \var{operationModes}.\func{add}(\var{bestMode})\;
    }
  }
  return \var{operationModes}\;
  \caption{Initial solution creation}
  \label{algo:initialSolve}
\end{algorithm}

  There are a few things to note about Algorithm~\ref{algo:initialSolve}.
  First, the parameter \texttt{prevMode} given to the calls for solving the models \MIPsf in line~\ref{line:init_callMIPsf} and \MIPs in line~\ref{line:init_callMIPs} is the operation mode of the previous time step, which we need to determine the operation mode change and compressor unit start variables \opModeChg{t} respectively \unitStart{u,t} for some unit $u\in\setCompressorUnitsI{a}\,a\in\setCompressorStations$ and $t\in\setTimestepsNoZero$ as explained in Section~\ref{sec:modelVariants}.
  In the call to \MIPs we furthermore give the parameter \texttt{validModes}, which replaces the set of valid network station operation modes \setNSModes.
  We also call the functions \texttt{modeAvailable} and \texttt{transitionsWork} in Algorithm~\ref{algo:initialSolve}.
  These refer to the operation mode unavailability and the transition time restriction introduced in Section~\ref{sec:gasNetworkStationDefinition}, where \texttt{modeAvailable} checks for a given operation mode $o$ and time $t$ if $o$ is available at $t$ and \texttt{transitionsWork} performs the checks described in Section~\ref{sec:gasNetworkStationDefinition} to test if a given operation mode sequence is valid regarding the corresponding transition times.
  In addition, Algorithm~\ref{algo:initialSolve} uses a function called \texttt{notSoonInfeasible}.
  Here, we check if choosing a new operation mode for time $t$ would result in an infeasibility at one of the subsequent time steps caused by a combination of operation mode unavailability and too long transition times.
  More specifically, we check if the new operation mode for time $t$ will become unavailable in one of the future time steps.
  If this is the case, we then test if there is enough time left to transition into another operation mode until then, also taking into account the time needed by the transition from the old operation mode at time $t-1$ to this new operation mode at time $t$.
  This look into the near future turned out to be necessary according to our computational experiments in order to avoid cases where Algorithm~\ref{algo:initialSolve} gets stuck in infeasible situations.

  Two other lines in the algorithm may require further explanation.
  In line~\ref{line:init_checkLastOpMode} we decide if the previous operation mode is good enough for the current time step by comparing its stationary objective function value against the cost of an operation mode change \costOpModes.
  If the objective function value is indeed smaller we know that the previous operation mode is the best option considering this individual time step given the chosen operation modes for the past time steps,
  since each other operation mode would at least have to pay the penalty of \costOpModes for changing the operation mode.
  As a last point to mention, Algorithm~\ref{algo:initialSolve} may abort without a feasible solution in line~\ref{line:abortInitialSolve}.
  This is no proof of infeasibility, since in theory cases are possible, where we abort although a feasible solution exists.
  However, in all of our test cases we have never aborted the algorithm at this point, also see our computational results in Section~\ref{sec:computational_results_results}.
  Furthermore, the combination of transition times and unavailable operation modes can lead to very hard to find feasible solutions, which makes the design of an algorithm performing reasonably fast on average but guaranteeing to find all feasible solutions a challenge.
  Therefore, we leave this problem open for future research.

\subsubsection{An improvement heuristic}
\label{sec:localImprovement}

  After we have found a feasible solution using Algorithm~\ref{algo:initialSolve}, we look for further improvements of it.
  Due to its design, Algorithm~\ref{algo:initialSolve} only considers individual time steps to decide which operation mode to choose for each time step.
  However, we can easily imagine a situation in which the operation mode $o_1$ found by \MIPs is best for time step $t$, but another operation mode $o_2$ is slightly better for all subsequent time steps and would have been the overall better choice at time $t$.
  We might even be able to avoid operation mode changes, if $o_1$ would become unavailable in the future, while $o_2$ has only slightly worse objective function values, but stays available.

\begin{algorithm}[ht]
  \DontPrintSemicolon
  \KwData{A sequence $S$ of valid operation modes over time}
  \KwResult{A valid sequence $S^*$ of operation modes with $\func{obj}(S^*)\leq\func{obj}(S)$}
  $\var{backwards} \gets True$\;
  \While{not having two iterations without improvements}{
    $\var{changeTimes} \gets$ list of times $t$ with $S[t-1]\neq S[t]$\;
    \lIf{\var{backwards}}{$\func{reverse}(\var{changeTimes})$}
    \For{$t\in \var{changeTimes}$\label{line:improve_changeTimeLoop}}{
      \lIf{$S[t-1]=S[t]$}{continue\label{line:improve_continue}}
      \eIf{\var{backwards}}{
        $\var{phaseToReplace} \gets$ phase ending at $S[t-1]$\;
      }
      {
        $\var{phaseToReplace} \gets$ phase starting from $S[t]$\;
      }
      $S^\text{best}, \var{bestImprovement} \gets$ \func{list}(), $0.0$\;
      \For{$\var{newMode} \in \func{convexCombination}(S[t-1], S[t])$}{
        $S^\text{new} \gets S.\func{replace}(\var{phaseToReplace}, \var{newMode})$\;
        $\var{improvement} \gets \func{obj}(S) - \func{obj}(S^\text{new})$\;
        \If{$\func{allModesAvailable}(S^\text{new})$
        \nonl\algoAnd $\func{transitionsWork}(S^\text{new})$
        \nonl\algoAnd $\var{improvement} > \var{bestImprovement}$}{
          $S^\text{best}, \var{bestImprovement} \gets S^\text{new}, \var{improvement}$\;
        }
      }
      \If{$\var{bestImprovement} >0$}{
        \tcp{Found improvement in this interation!}
        $S \gets S^\text{best}$\;
      }
    }
    $\var{backwards} \gets not \var{ backwards}$\;
  }
  return $S$\;
  \caption{Improvement heuristic}
  \label{algo:localImprovement}
\end{algorithm}

  To deal with these situations, we created for a given feasible solution, represented by a sequence of operation modes over time, the improvement heuristic stated as Algorithm~\ref{algo:localImprovement}.
  Here the idea is, to identify all sequences of identical operation modes over time in the solution.
  We call these sequences stable phases or just \emph{phases} of a feasible solution.
  Obviously, the switch from one phase to the subsequent one happens if the operation mode changes to a new mode.
  For each of these phases we then check if we can replace the operation mode of the whole phase with a similar one being more beneficial in terms of the objective function value.

  We obtain these similar network station operation modes from the call of the function \texttt{convexCombination}, which is the key feature of Algorithm~\ref{algo:localImprovement}.
  To define it, we use the the function $M(o,a)$ returning the mode or active configuration of a valve or compressor station $a$ in operation mode $o$, see Section~\ref{sec:gasNetworkStationDefinition}.
  Furthermore, we denote by $\setCompressorUnits(x)$ the compressor units used in mode or configuration $x\in\{\text{by}, \text{cl}\} \cup \setCompressorConfigurations{a}$ for some compressor station $a\in\setCompressorStations$, where $\setCompressorUnits(\text{by}) = \setCompressorUnits(\text{cl}) = \emptyset$.
  Then we first define the function \texttt{convexCombinationCS} on a tuple $(x,y)$ with $x,y\in\{\text{by}, \text{cl}\} \cup \setCompressorConfigurations{a}$ as
  \begin{align*}
    \text{\texttt{convexCombinationCS}}(x,y) &:= \{ x,y \} \cup \\
    \Big\{ c\in\setCompressorConfigurations{a} \quad | \quad & \forall u\in \setCompressorUnits(x)\cap\setCompressorUnits(y): u\in\setCompressorUnits(c)\\
    \land\, & \forall u\in \setCompressorUnits(c): u \in\setCompressorUnits(x)\cup\setCompressorUnits(y) \Big\}.
  \end{align*}
  Note, that \texttt{convexCombinationCS}$(x,y) \subseteq \{\text{by}, \text{cl}\} \cup \setCompressorConfigurations{a}$ holds.
  Then we are ready to define \texttt{convexCombination} on a tuple $(o_1,o_2)$ of operation modes as
  \begin{align*}
    \text{\texttt{convexCombination}}(o_1,o_2) :=& \Big\{ o\in\setNSModes \,| \\
    \big( \forall a\in\setValves: M(o,a) =& M(o_1,a) \quad \lor \quad \forall a\in\setValves: M(o,a) = M(o_2,a) \big)\\
    \land\quad \forall a\in\setCompressorStations: M(o,a) \in& \text{\texttt{convexCombinationCS}}\big(M(o_1,a), M(o_2,a)\big) \Big\}.
  \end{align*}
  Note here, that while we allow a compressor station $a$ to have a configuration using a compressor unit set ``in between'' the used compressor unit set of the configurations used in $o_1$ and $o_2$ for $a$, we only allow the exact valve mode combination used in $o_1$ or the one used in $o_2$.
  The reason for this is, that a valve mode combination enables a very specify set of paths through each network station, and it is very unlikely that a valve mode set obtained from combining the modes used in the two given operation modes yields operation modes, which are able to handle the same demand situation.
  Since, each of the modes obtained by calling \texttt{convexCombination} is tested in Algorithm~\ref{algo:localImprovement}, we hereby restrict the result set to the most promising candidates.

  Apart from calling \texttt{convexCombination}, Algorithm~\ref{algo:localImprovement} uses the two functions \texttt{transitionsWork}, which works in the same way as described for Algorithm~\ref{algo:initialSolve} above and \texttt{allModesAvailable}, which is similar to \texttt{modeAvailable} from Algorithm~\ref{algo:initialSolve}, but instead of checking the availability of a given operation mode $o$ for time $t$ checks the availability of a whole sequence of operation modes at the times corresponding to the position in the sequence.
  Furthermore, we evaluated the objective function value of a sequence of operation modes using the function \texttt{obj} by successively calling \MIPsf for each operation mode and time corresponding to its position in the sequence.
  If one of the models turns out to be infeasible, the returned objective function value will be infinity.

  Finally, we note that we decided to start the algorithm in the backwards oriented mode.
  The reason for this decision is that the initial solution is obtained by Algorithm~\ref{algo:initialSolve} which was operated in a forward direction.
  Furthermore, we highlight that Algorithm~\ref{algo:localImprovement} has the potential to reduce the total number of needed operation mode changes, since the two original operation modes are always part of the result of \texttt{convexCombination}.
  In addition, it is possible that an operation mode change from the loop of line~\ref{line:improve_changeTimeLoop} has already been removed in the previous iteration by replacing one of the involved operation modes with the other one, which makes the check in line~\ref{line:improve_continue} necessary.

\subsection{Transient solution smoothing}
\label{sec:smoothing}

  As a final step of our specialized network station algorithm, we solve the transient model variant \MIPf with fixed network station operation modes and flow directions.
  We obtain both for each time step from the stationary model solutions created in the previous steps of the algorithm.
  When comparing pressure and flow values at single points of network over time, we expect the transient solution states to be more similar in general and in case of changing conditions to be more smooth compared to the series of stationary solution states.
  This is due to the missing penalty of operation point changes in the independent stationary models, which may result in considerably different solution states.
  This difference may for example occur in the pressure level of nodes inside the station, even if the demand situation as well as the determined operation mode and flow direction are the same.

\begin{algorithm}[ht]
  \DontPrintSemicolon
  \KwData{A sequence $S$ of tuples of operation modes and flow directions over time as well as a time horizon size h}
  \KwResult{A set of transient solution states for all $t\in\setTimesteps$}
  \If{$\var{h} \geq |\setTimesteps|$}{
    \tcp{overall time horizon covered by smoothing time horizon}
    return $\MIPf(S, initialState = \var{initialState})$\;
  }
  $\var{solutionStates} \gets \func{list}()$\;
  $\var{solutionStates}.\func{add}(\var{initialState})$\;
  $\var{currTime} \gets 0$\;
  \While{$\var{currTime} + h \leq |\setTimesteps|$}{
    $S^\text{h} \gets S.\func{slice}(\var{currTime}, \var{currTime} + \var{h})$\;
    $\var{thisTimeStates} \gets \MIPf (S^\text{h}, initialState = \var{solutionStates}.\func{last}())$\;
    \eIf{$\var{currTime} + h = |\setTimesteps|$}{
      $\var{solutionStates}.\func{addAll}(\var{thisTimeStates})$\;
    }{
      $\var{solutionStates}.\func{add}(\var{thisTimeStates}.\func{first}())$\;
    }
    $\var{currTime} \gets \var{currTime} + 1$\;
  }
  return \var{solutionStates};
  \caption{Transient smoothing}
  \label{algo:smoothing}
\end{algorithm}

  In our computational experiments we observed that even though most of the binary decision variables of \MIP are fixed in \MIPf, only a limited number of time steps can be solved for large network stations.
  Therefore, we use a \emph{rolling horizon} approach to solve \MIPf, which is described in Algorithm~\ref{algo:smoothing}.
  Here, we specify a small fixed time horizon size $h$, which represents the number of time steps to solve in model \MIPf including the time step for the given initial state.
  We then solve a series of models \MIPf, while always fixing the earliest time step and shifting the time horizon by 1 in each iteration.
  In the function call to solve \MIPf in Algorithm~\ref{algo:smoothing}, we give the subsequence of operation modes and flow directions, corresponding to and also encoding the current time horizon to solve.
  Furthermore, we specify the state to use as the fixed initial state.

  The main benefit of this method is, that increasing the size $|\setTimestepsAll|$ of the overall time horizon only increases the number of equally sized and therefore similarly complex MIP models to solve rather then increasing the complexity of the model, which may lead to an exponential increase in runtime.

\section{Computational experiments}
\label{sec:computational_results}
  In order to verify the competitiveness of Algorithm~\ref{algo:completeStationModel} in terms of both solution quality as well as execution time, we evaluated it to a large number of test instances.
  These instances represent network stations in the network of our project partner Open Grid Europe(OGE), for which we generated scenario values based on historic real-world situations in the network.

  Throughout this chapter, we will state all flow values as volumetric flow under normal conditions in the unit $1000\text{m}^3/\text{h}$.
  To create these values from mass flow values in kg/s we use the normal density \nDens given in kg/m$^3$, which is determined by the gas mixture and assumed to be constant all over the network.
  The linear transformation formula then reads as $[1000\text{m}^3/\text{h}] = \frac{3600}{1000\nDens} [\text{kg}/\text{s}]$.

\subsection{Instances}
  \label{sec:instanceDescription}
  We consider 7 different network stations from the network of OGE with different sizes and properties.
  An overview can be found in Table~\ref{tab:stationStatistics}.
  For each of these stations, we have a set of 159 instances.
  These have been created by solving a macroscopic gas flow model on an aggregated complete network containing simplified versions all the network stations.
  A full description of the model can be found in~\cite{HopHenLenKoc2019}.

  \begin{table}[ht]
    \centering
    \begin{tabular}{crrrrrrrr}
      Name & $|\setVertices|$ & $|\setArcs|$ & $\frac{\sum\limits_{a\in\setPipes}\lenI{a}}{|\setPipes|}$ &
      $|\setCompressorConfigurations{a}| \enskip \forall a\in\setCompressorStations$ &
      $|\setNSModes|$ & $|\setFlowDirectionsNavi|$ & $|\setBoundaryNodesExitP|$ & $|\setFlowDirectionConditions|$\\
      \midrule
      A              &  14 &  11 & 0.001 km &             16 &   34 &  3 & 1 & 0\\
      B              &  11 &  12 & 0.001 km &             18 &   23 &  2 & 1 & 0\\
      C              &  27 &  34 & 0.012 km &              2 &   13 &  4 & 2 & 0\\
      D              &  25 &  31 & 0.404 km &           2, 6 &   92 &  6 & 2 & 0\\
      E              &  48 &  67 & 0.308 km &           3, 5 &   82 & 12 & 1 & 4\\
      F              &  51 &  66 & 0.024 km &    2, 3, 6, 12 & 2836 &  3 & 3 & 0\\
      G$_\text{min}$ & 118 & 148 & 0.079 km & 1, 1, 2, 7, 20 & 1267 & 15 & 2 & 0 \\
      G$_\text{max}$ & 120 & 150 & 0.079 km & 1, 1, 2, 7, 20 & 1285 & 20 & 2 & 0
    \end{tabular}
    \caption{Overview of different properties of the 7 network stations A to G.
    For station G the topology changed during the considered time period of 91 days. Hence, we denoted by G$_\text{min}$ the minimal values of all the quantities for station G and by G$_\text{max}$ the corresponding maximum values.}
    \label{tab:stationStatistics}
  \end{table}

  As initial states we used state values that actually occurred for the network of our project partner.
  Analogously we are given measured pressure and flow values for the next 12 hours after the initial state time at the boundary nodes of the network containing all the 7 network stations.
  These instances are distributed over a period of 91 days.
  For each day we have two instances, whose initial state times have a difference of 12 hours.
  Inside this data period, we unfortunately are missing instances for 8 days as well as having only one instance for another 7 days due to technical problems during the data creation at our project partner.
  Therefore, our final instance set consists of $2\cdot 91 - 2\cdot8 - 7 = 159$ instances.

  The solutions of the macroscopic gas flow model for each of these instances yield pressure and inflow values at the boundaries of each network station.
  These values represent the scenario values for each instance of the test set of this paper.
  For the macroscopic gas flow model, the time horizon of 12 hours has been split into 4 periods of 15 minutes followed by 11 periods of 60 minutes, resulting in 15 future time steps in total.

  For our analysis, we created the following four different partitions of the 12 hour time horizon, where the horizon is built from left to right per row.
  \begin{align*}
    12 \text{ time steps:} \qquad &  \makebox[\widthof{48}][r]{4}\times \makebox[\widthof{7.5}][r]{15}\,\text{min, }  \makebox[\widthof{18}][r]{5}\times 60\,\text{min, } 3\times 120\,\text{min} \\
    24 \text{ time steps:} \qquad & \makebox[\widthof{48}][r]{4}\times \makebox[\widthof{7.5}][r]{15}\,\text{min, } 18\times 30\,\text{min, } 2\times \makebox[\widthof{120}][r]{60}\,\text{min} \\
    48 \text{ time steps:} \qquad & 48 \times \makebox[\widthof{7.5}][r]{15}\,\text{min} \\
    96 \text{ time steps:} \qquad & 96 \times 7.5\,\text{min}
  \end{align*}
  Our partitioning of the time horizon can have additional or missing time points compared with the 15 time steps used in the macroscopic model.
  In order to create scenarios values that fit the time horizon we interpolate using the original values.

  To conclude, we used the following set of weights for the objective function for all instances:
  \begingroup
  \addtolength{\jot}{0.5em}
  \begin{align*}
      \costSlackP        &= \frac{1000.0}{\text{bar}\cdot\text{h}} &
      \costOpModes       &=       1000.0 &
      \costRgPLChg   &= \frac{  10.0}{\text{bar}} &
      \costCsPLChg   &= \frac{  10.0}{\text{bar}} \\
      \costSlackQ        &= \frac{ 100.0}{1000\,\text{m}^3} &
      \costUnitStarts    &=       1200.0 &
      \costRgPRChg   &= \frac{  10.0}{\text{bar}} &
      \costCsPRChg   &= \frac{  10.0}{\text{bar}} \\
      &&
      \costRgModeChg &=         50.0 &
      \costRgQChg    &= \frac{   1.0}{1000\,\frac{\text{m}^3}{\text{h}}} &
      \costCsQChg    &= \frac{   1.0}{1000\,\frac{\text{m}^3}{\text{h}}}
  \end{align*}
  \endgroup
  Note that the unit of \costSlackQ is deduced by $\frac{1}{1000\,\frac{\text{m}^3}{\text{h}}\cdot\text{h}} = \frac{1}{1000\,\text{m}^3}$.

\subsection{Computational setup}

  We performed our computations on a cluster using 4 cores and 16\,GB of RAM of a machine being composed of two \emph{Intel Xeon CPU E5-2670 v2} running at 2.50\,GHz.
  As a solver for the underlying MIP problems we used \emph{Gurobi} in version 8.1.0~\cite{Gur2018}, which we accessed via the Pyomo modeling language~\cite{HarLaiWatWoo2017}\cite{HarWatWoo2011}.
  Since the corresponding MIP models are numerically challenging we used the solver with maximal  \emph{NumericFocus} parameter.
  In addition, we specified in Table~\ref{tab:MIPsettings} the optimality conditions and maximum run times for each solved model variant.
  \begin{table}[ht]
  \centering
  \begin{tabular}{lrrr}
      Variant       & Rel. Gap  & Abs. Gap  & Time limit \\
      \MIPs, \MIPsf & 1E-4      & 1E-2      & 10h \\
      \MIPf         & 5E-3      & 1E-2      & 60s
  \end{tabular}
  \caption{Optimality conditions and maximum run times we use for the single model variants. Note that the 10 hours serve as a representation for $\infty$, i.e., we have chosen the time limit high enough to always solve the models to optimality.}
  \label{tab:MIPsettings}
  \end{table}

  Finally, we specify the rolling horizon parameter $h$ to be 4, so we are always solving the smoothing for 4 future time steps.
  We found in our experiments, that this number represents a good trade-off between efficient model solving speed and foresightedness of the obtained solution.

\subsection{Results}
\label{sec:computational_results_results}

  As described in Section~\ref{sec:instanceDescription} above, we tested Algorithm~\ref{algo:completeStationModel} on $159 \text{ start times} \times 7 \text{ network stations} \times 4 \text{ time horizons } = 4452$ instances.
  For all of these instances, it found a feasible solution.
  In particular, we always found a feasible solution in the initial solution creation in Algorithm~\ref{algo:initialSolve} described in Section~\ref{sec:initialSolution} although this is not guaranteed by the algorithm itself.

  \begin{figure}[ht]
      \centering
      \includegraphics[trim={1.0cm 0.0cm 2.0cm 1.0cm},clip,width=0.8\linewidth]{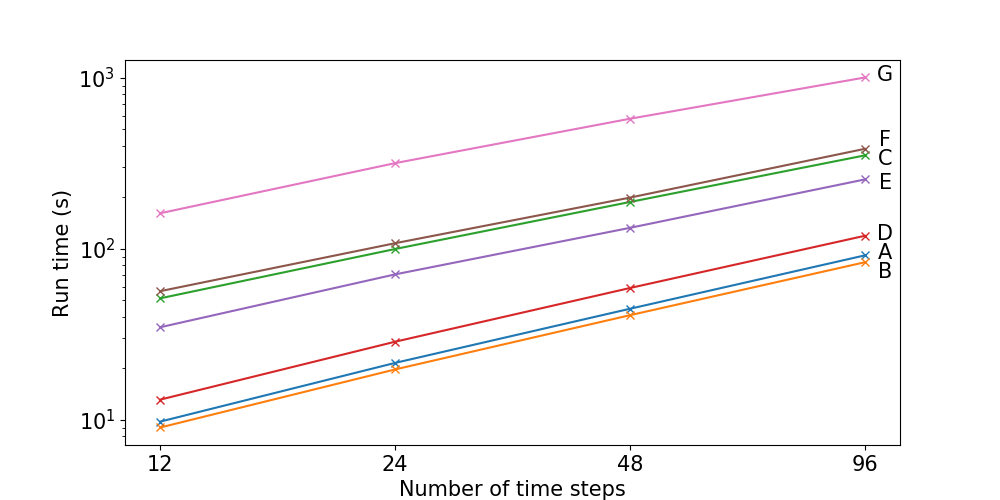}
      \caption{Average run times for Algorithm~\ref{algo:completeStationModel} sorted by number of time steps and network station.}
      \label{fig:logRunTimes}
  \end{figure}

  In Figure~\ref{fig:logRunTimes} one can find an overview of the average run times of Algorithm~\ref{algo:completeStationModel} sorted by station and number of solved time steps.
  As one can see, the ranking of the single stations regarding run time as well as the ratio between the single stations run times is stable over the different time horizon partitions.
  Furthermore, the run time per station corresponds roughly to the complexity deduced from the station statistics in the sense that more nodes, arcs and operation modes make the overall problem more complex to solve and therefore increase the run time.
  The exception to this reasoning is station C, which has a rather high average run time although its node and arc counts are comparable to station D and the number of operation modes is even the smallest of all seven stations.
  We assume that this is due to the stations topology, in which the regulators seem to be arranged in a way which is hard to solve for the final smoothing step of the algorithm.
  Another noticeable fact is the almost perfectly proportional average run time increase with increasing number of time steps.
  As mentioned in Section~\ref{sec:smoothing} this property was an explicit goal when designing the algorithm.

\begin{figure}[ht]
      \centering
      \includegraphics[trim={1.0cm 0.8cm 2.0cm 1.5cm},clip,width=0.8\linewidth]{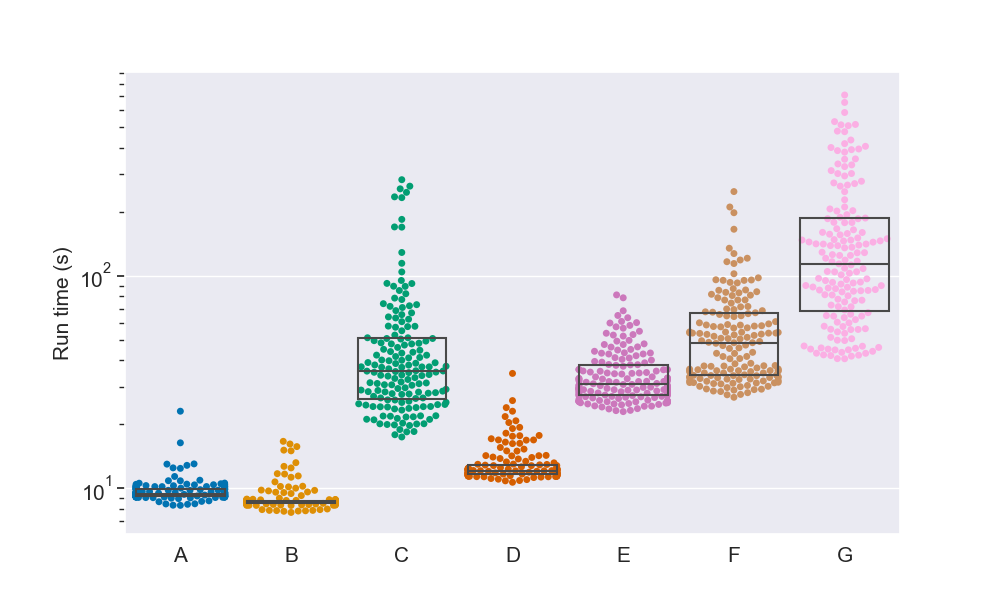}
      \caption{Run times of instances with 12 time steps, sorted by station. One dot represents one instance, however dots may overlap when too many have similar run time. In addition we draw a typical whiskerless box plot. Therefore the lines represent the 25th percentile, the median and the 75th percentile.}
      \label{fig:spreadRunTimes}
  \end{figure}

  In addition to the average run times, we also depict the distribution of run times sorted by station in Figure~\ref{fig:spreadRunTimes} for instances solved for 12 time steps.
  The distribution of the other time steps are quite similar in appearance and can be found in the Appendix as Figures~\ref{fig:spreadRunTimes24}, \ref{fig:spreadRunTimes48}, and \ref{fig:spreadRunTimes96} respectively.

  For the fast solving stations A, B, and D, the run times of the instances are all very similar in that the middle box is barely visible.
  Station C is again outstanding by having a rather large spread of run times only matched by the biggest and most complex station G.
  However, in general the difference between the median and maximum run time are lower than an order of magnitude for all of the stations.
  That means that for our instance set even the extreme cases are still in reach and have a somewhat similar run time to the rest of the instances, making Algorithm~\ref{algo:completeStationModel} well suited for strict time limit situations, which we would face in a production environment.

  \begin{figure}[tb]
    \centering

    \begin{subfigure}[c]{0.48\textwidth}
    \includegraphics[trim={1.5cm 0.5cm 1.5cm 0.0cm},clip,width=\textwidth]{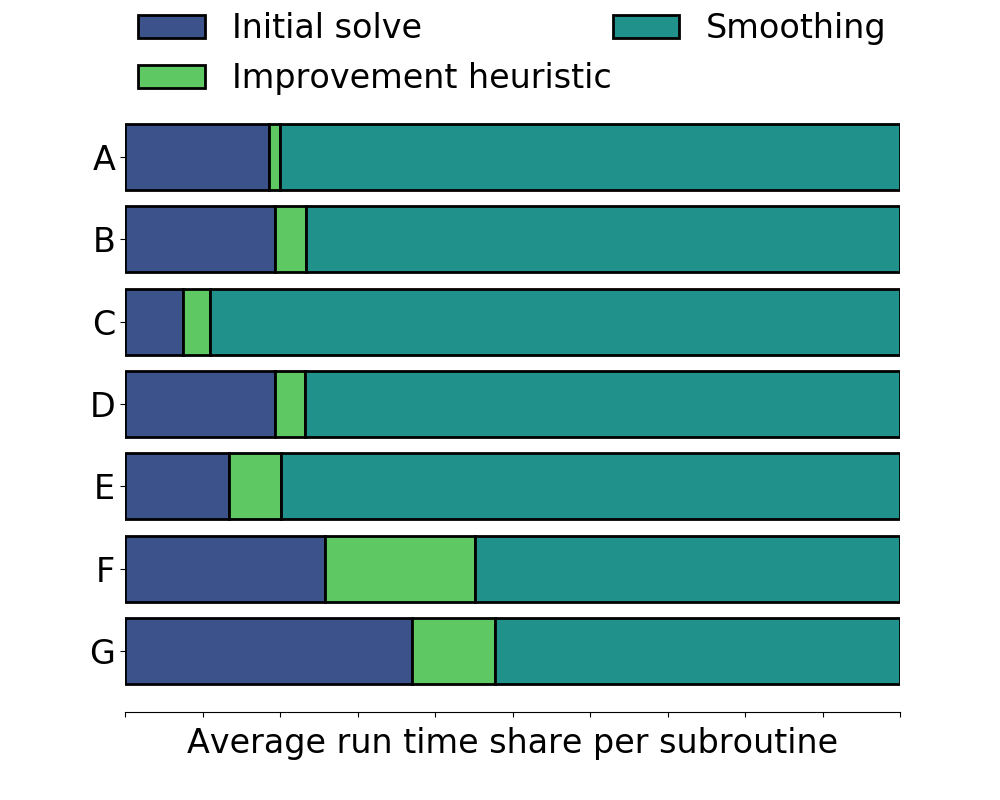}
    \subcaption{12 Timesteps}
    \label{fig:runTimeShare12}
    \end{subfigure}
    \begin{subfigure}[c]{0.48\textwidth}
    \includegraphics[trim={1.5cm 0.5cm 1.5cm 0.0cm},clip,width=\textwidth]{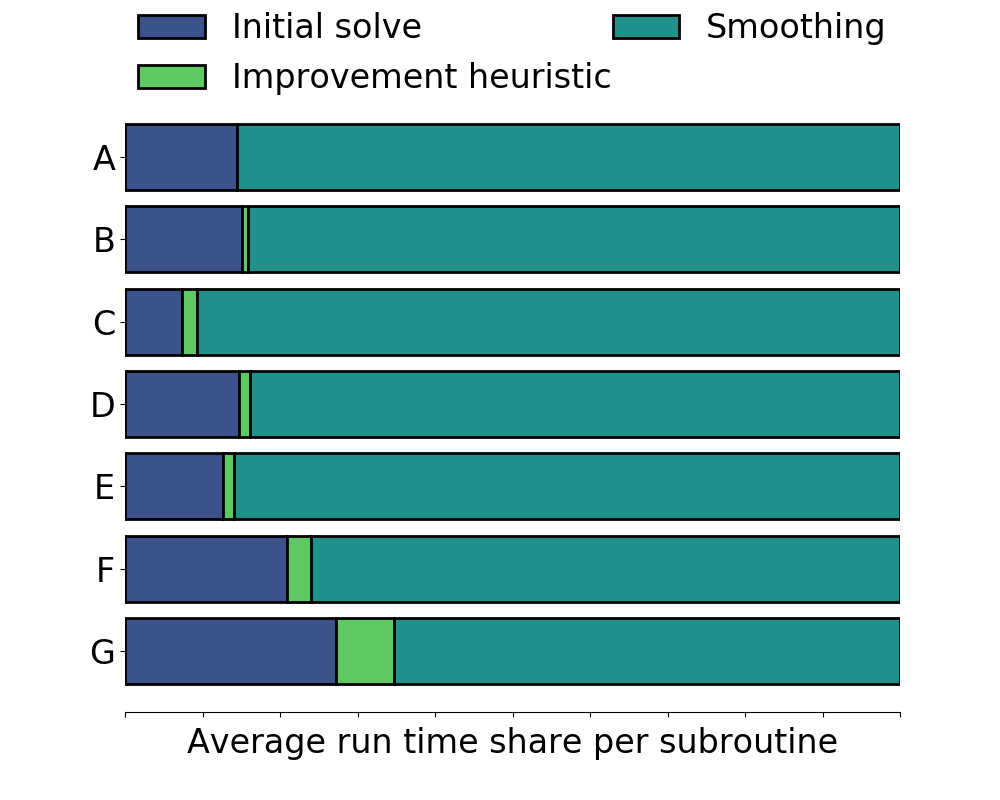}
    \subcaption{96 Timesteps}
    \label{fig:runTimeShare96}
    \end{subfigure}
    \caption{Portion of run time spend in the single subroutines of Algorithm~\ref{algo:completeStationModel} displayed for each of the single stations. Figure~\ref{fig:runTimeShare12} represents the results for 12 time steps, Figure~\ref{fig:runTimeShare96} the results for 96 time steps.}
    \label{fig:runTimeShares1296}
  \end{figure}

  As the last run time related graphic, Figure~\ref{fig:runTimeShares1296} shows the average portion of run time spent in the single subroutines of Algorithm~\ref{algo:completeStationModel} for 12 time steps as well as for 96 time steps.
  It states that the transient smoothing dominates the run time already for 12 time steps.
  With an increasing number of time steps, the influence even increases, while the improvement heuristics impact is decreasing.
  This can also be verified through the 24 and 48 time step instances displayed in Figure~\ref{fig:runTimeShares2448} in the Appendix.
  The observed shares fit to the structure of the single subroutines from Section~\ref{sec:stationModelAlgorithm}.
  While the initial solve as well as the smoothing have to solve more instances of the challenging model variants \MIPs respectively \MIPf with each additional time step, the improvement heuristic iterates over pairs of subsequent time steps with different operation modes.
  Since the differently sized time horizon partitions all cover the same 12 hours, the number of total operation mode changes in the result of the initial solve is very unlikely to increase much with increasing time granularity.
  This is also due to the initial solve algorithm itself, which tries to use the same operation mode as long as possible before changing to a different one.

  Apart from the run time we are also interested in the quality of the solution of Algorithm~\ref{algo:completeStationModel}.
  Therefore, we solved the complete model \MIP directly on the smaller instances A to E to obtain a valid lower bound for them.
  Note that \MIP itself is already only a relaxation of the problem solved by Algorithm~\ref{algo:completeStationModel}, since the transition time restrictions for operation modes mentioned in Section~\ref{sec:gasNetworkStationDefinition} are not included in this model.
  The obtained lower bound can be used to give an upper bound on the relative difference of the solution found by Algorithm~\ref{algo:completeStationModel} to the optimal solution.
  We call this difference \emph{gap} and define it as $\frac{\text{obj}(\text{solution}) - \text{lowerBound}}{\text{obj}(\text{solution})}$.
  Since the lower bound cannot be negative, the gap is always smaller than or equal to 100\%.

  We solved model \MIP on stations A to E for 12 time steps using the same computational setup as described above with a time limit of 10 hours.
  Furthermore, we gave the solution obtained from Algorithm~\ref{algo:completeStationModel} as a starting solution.
  The solver finished all 159 instances of the stations A and B within the time limit.
  For the stations C, D and E we hit the time limit for 47, 8 and 18 instances respectively, resulting in potentially sub-optimal lower bounds in these cases.
  For the stations F and G more than half of the instances did not finish in time, which is why we did not include these in our evaluation.

  \begin{figure}[ht]
      \centering
      \includegraphics[width=0.7\linewidth]{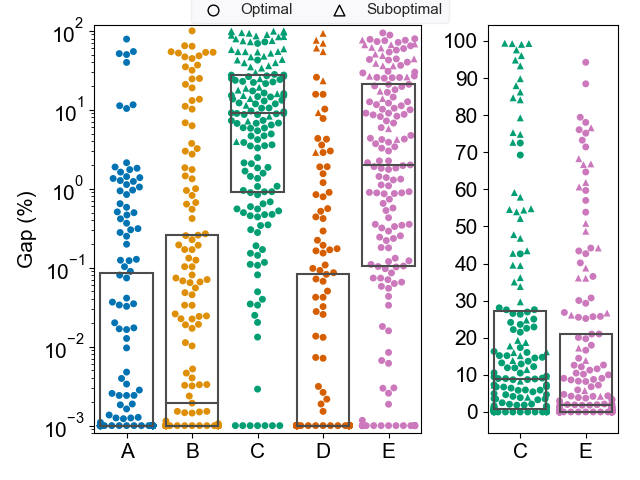}
      \caption{Gap between the solution of Algorithm~\ref{algo:completeStationModel} and a lower bound obtained from solving \MIP directly.
      We considered the stations A to E and solved them using 12 time steps.
      One dot represents one instance, however dots may overlap when too many have similar run time.
      The instances with potentially suboptimal lower bounds due to hitting the time limit are marked by triangles instead of dots.
      In addition we draw a typical whiskerless box plot.
      Therefore the lines represent the 25th percentile, the median and the 75th percentile.
      The left graphic uses logarithmic scale for the gap and the right one displays the same data for stations C and E on linear scale.
      For all instances in which the objective function value and the lower bound are smaller than 0.1 we defined the gap to be zero.
      Furthermore we plotted all values below $10^{-3}$ as $10^{-3}$ for the logarithmic scale.}
      \label{fig:stationModelOptGap}
  \end{figure}

  The results can be found in Figure~\ref{fig:stationModelOptGap}, where we plotted the gap for all 5 network stations on a logarithmic scale and then plotted the same data for stations C and E on a linear scale to better depict their distributions.
  Note that the instances which have a potentially suboptimal lower bound are marked with triangles instead of dots.
  As one can observe, we have good results for all of the network stations with more than half of the instances for each station having a solution with at most 10\% difference to the lower bound.
  For the best three stations A, B, and D we can even state that at least 75\% the of values have less than 1\% difference to the optimal solution.
  The picture is more diverse for the other two stations, in which the 75\% percentile is between 20\% and 30\% gap.
  However, for those stations we have the highest value of instances for which \MIP did not solve the problem in time.
  Therefore, there is the potential for further improvement of the gap by increasing the lower bound here.
  For all stations C, D and E with instances having potentially suboptimal lower bounds, the corresponding instances are among those having high gap values.
  The picture is extreme for Station C, where all but 3 instances above the 75\% percentile had suboptimal lower bounds, making it very likely that one can further decrease the gap value by increasing the run time limit of \MIP directly.

  Altogether, the results show that Algorithm~\ref{algo:completeStationModel} is able to find good solutions in a short time and is therefore able to solve even large instances of the presented gas transportation problem, which would be far out of reach when trying to solve the original MIP model formulation \MIP directly.

\section{Conclusion}
\label{sec:conclusion}

  We presented the transient gas network transportation problem on so-called network stations, which represent the intersection points of major transportation pipelines and contain the majority of active elements to control the network.
  For this problem we introduced a mixed-integer programming model including a complex model for compressor stations as well as additional variables and constraints for the network station itself.
  For the pipes, we found that due to their shortness they have less overall impact in network stations, which enabled us to use a linear description for them.
  Using state-of-the art solvers the MIP Model is not tractable for large network stations.
  Therefore we developed a specialized algorithm to solve the problem.
  Here we again use the fact that network stations contain only short pipes and therefore negligible linepack, causing the decision making to mostly depend on the current demand situation at the boundary of the station.
  Therefore, we determined the operation modes and gas flow directions of the network station based on multiple solves of a stationary version of the MIP.
  By using this approach we are also able to satisfy transition time constraints for the operation modes, which have been excluded from the original MIP formulation.
  In order to obtain a feasible solution for the problem, we finally solved another variant of the MIP with fixed network station operation modes and flow directions in a rolling horizon fashion.

  To verify the competitiveness of our algorithm regarding both run time and solution quality, we did tests on 159 different scenarios based on past flow situations in 7 real world network stations provided by our project partner.
  Our algorithm is able to compute feasible solutions for all of the presented instances very fast.
  Even for the biggest of the presented network stations consisting of more than 100 nodes and arcs as well as over 1000 different network station operation modes, our algorithm terminates on average in under 20 minutes for each of the stations.
  By running experiments using 4 different types of granularity, we found that the run time increases proportional to the number of time steps used, indicating that our algorithm scales very well.
  In terms of quality of the solution, we tested the results against a lower bound obtained from solving the original MIP formulation for 12 time steps on all except the biggest 2 network stations.
  For the worst of the 5 remaining stations, more than half of the instances have a solution with at most 10\% difference to the lower bound and for the best 3 stations, we find near optimal solution with less than 1\% difference for more than 75\% of the instances.
  Altogether we were able to show, that our algorithm can reliably find good solutions to the problem in a short amount of time.

  There are a lot of different possibilities to continue this research.
  To increase the model accuracy, the approximative linearization of the friction in pipes and resistors as well as the maximum power bound of compressor units could be replaced by their original non-linear versions.
  This would turn the model into a MINLP and thereby increase its complexity.
  From a theoretical point of view, extending Algorithm~\ref{algo:initialSolve} such that it has the guarantee to always find existing feasible solutions would greatly improve the overall robustness.
  Finally, real-world gas network operation is a complicated business with a never-ending list of special elements and extra constraints, which can still be added to our model.
  As examples we name ramp-up and cool down times for compressor units as well as target value based control of regulators and compressor stations.

\section*{Acknowledgements}
The work for this article has been conducted in the Research Campus MODAL funded by the German Federal Ministry of Education and Research (BMBF) (fund number 05M14ZAM).

\bibliographystyle{abbrv}
\bibliography{bibliography.bib}

\newpage
\appendix
\renewcommand\thefigure{\thesection.\arabic{figure}}
\setcounter{figure}{0}
\renewcommand\thetable{\thesection.\arabic{table}}
\setcounter{table}{0}
\section{Appendix}

\begin{table}[htp]
\centering
\newcommand{\realPlus}{\hfill\makebox[\widthof{$\in\{0,1\}$}][l]{$\in\setRealsNonNeg$}}
\newcommand{\realAll}{\hfill\makebox[\widthof{$\in\{0,1\}$}][l]{$\in\setReals$}}
\newcommand{\binary}{\hfill$\in\{0,1\}$}

\begin{tabular}{lll}
Variable & Meaning & Unit \\
\hline
\pressI{v,t}             \realPlus & pressure at node $v\in\setVertices$ & bar \\
\mFlowI{a,t}              \realAll  & flow over arc $a\in\setArcs\setminus\setPipes$ & kg/s \\
\mFlowI{v,a,t}          \realAll  & flow into or out of pipe $a\in\setPipes$     & kg/s \\
\inflowI{v,t}          \realAll  & inflow into boundary node $v\in\setBoundaryNodes$ & kg/s \\
\cfg{c,a,t} \binary   & selection of configuration $c\in\setCompressorConfigurations{a}$ for $a\in\setCompressorStations$ & 1 \\
\modeBy{a,t}           \binary   & selection of bypass mode for $a\in\setRegulators\cup\setCompressorStations$ & 1 \\
\modeCl{a,t}           \binary   & selection of closed mode for $a\in\setRegulators\cup\setCompressorStations$ & 1 \\
\modeAc{a,t}           \binary   & selection of active mode for $a\in\setRegulators$ & 1 \\
\modeOp{a,t}             \binary   & selection of open/closed mode for $a\in\setValves$ & 1 \\
\pressBy{a,t}       \realPlus & pressure in bypass mode for $a\in\setCompressorStations$ & bar \\
\mFlowBy{a,t}           \realAll  & flow in bypass mode for $a\in\setCompressorStations$ & kg/s \\
\pressClL{c,a,t}   \realPlus & incoming pressure in closed mode for $a\in\setCompressorStations$ & bar \\
\pressClR{c,a,t}   \realPlus & outgoing pressure in closed mode for $a\in\setCompressorStations$ & bar \\
\pressCfgL{c,a,t}   \realPlus & incoming pressure in configuration $c\in\setCompressorConfigurations{a}$ of $a\in\setCompressorStations$ & bar \\
\pressCfgR{c,a,t}   \realPlus & outgoing pressure in configuration $c\in\setCompressorConfigurations{a}$ of $a\in\setCompressorStations$ & bar \\
\mFlowCfg{c,a,t}        \realAll  & flow in configuration $c\in\setCompressorConfigurations{a}$ of $a\in\setCompressorStations$ & kg/s \\
\opMode{o,t}               \binary   & selection of mode $o\in\setNSModes$ & 1 \\
\flowDir{f,t}            \binary   & selection of flow direction $\flowDirection\in\setFlowDirectionsNavi$ & 1 \\
\slackPPos{v,t}     \realPlus & positive pressure slack for boundary node $v\in\setBoundaryNodes$ & bar \\
\slackPNeg{v,t}     \realPlus & negative pressure slack for boundary node $v\in\setBoundaryNodes$ & bar \\
\slackQPos{v,t}         \realPlus & positive flow slack for boundary node $v\in\setBoundaryNodes$ & kg/s \\
\slackQNeg{v,t}         \realPlus & negative flow slack for boundary node $v\in\setBoundaryNodes$ & kg/s \\
\opModeChg{t}            \binary   & operation mode change & 1 \\
\regModeChg{a,t}        \binary   & mode change for regulator $a\in\setRegulators$ & 1 \\
\unitStart{u,t}            \binary   & start of compressor unit $u\in\setCompressorUnitsI{a}$ for $a\in\setCompressorStations$ & 1 \\
\rgPLChg{a,t}           \realPlus & change of incoming pressure of active $a\in\setRegulators$ & bar \\
\rgPRChg{a,t}           \realPlus & change of outgoing pressure of active $a\in\setRegulators$ & bar \\
\rgQChg{a,t}            \realPlus & change of flow over of active $a\in\setRegulators$ & kg/s \\
\csPLChg{a,t}           \realPlus & change of incoming pressure of active $a\in\setCompressorStations$ & bar \\
\csPRChg{a,t}           \realPlus & change of outgoing pressure of active $a\in\setCompressorStations$ & bar \\
\csQChg{a,t}            \realPlus & change of flow over of active $a\in\setCompressorStations$ & kg/s
\end{tabular}
\caption{List of all used variables, specifying their domain, meaning and unit. All variables are defined for $t\in\setTimestepsNoZero$, not taking into account the given values for the initial state, see Section~\ref{sec:scenario_initState}. Note that $1\,\text{bar}=10^5\,\text{Pa}$.}
\label{tab:Variables}
\end{table}

  \begin{figure}[htp]
      \centering
      \includegraphics[trim={1.0cm 0.8cm 2.0cm 1.5cm},clip,width=0.80\linewidth]{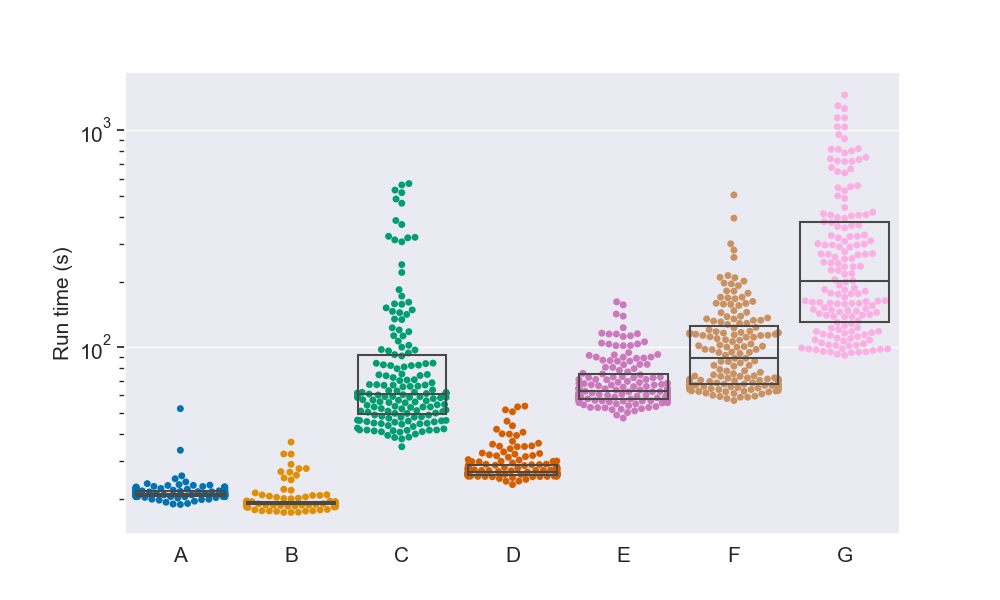}
      \caption{Run times of instances with 24 time steps, sorted by station. One dot represents one instance, however dots may overlap when too many have similar run time. In addition we draw a typical whiskerless box plot. Therefore the lines represent the 25th percentile, the median and the 75th percentile.}
      \label{fig:spreadRunTimes24}
  \end{figure}

  \begin{figure}[htp]
      \centering
      \includegraphics[trim={1.0cm 0.8cm 2.0cm 1.5cm},clip,width=0.80\linewidth]{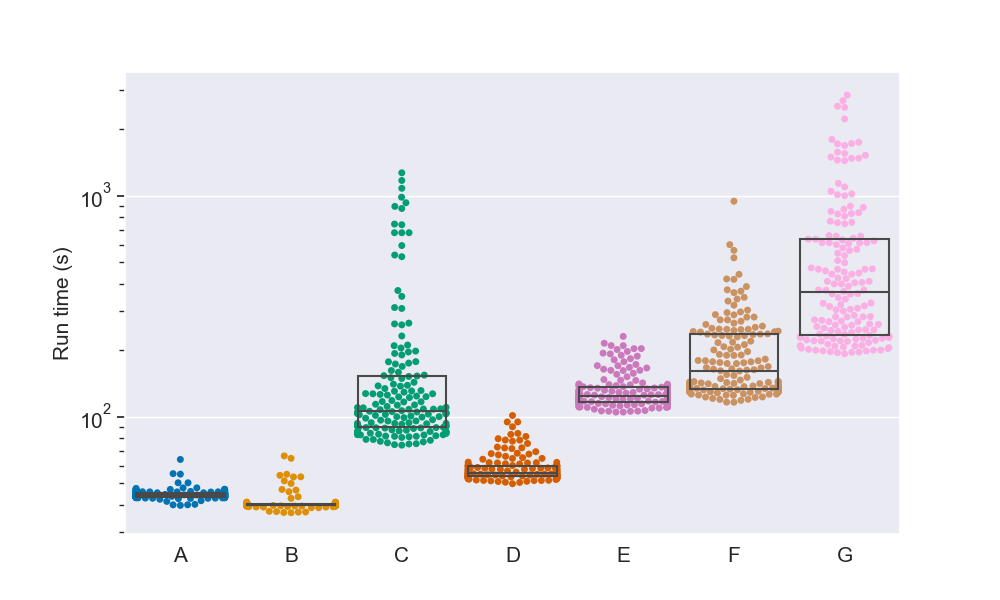}
      \caption{Run times of instances with 48 time steps, sorted by station. One dot represents one instance, however dots may overlap when too many have similar run time. In addition we draw a typical whiskerless box plot. Therefore the lines represent the 25th percentile, the median and the 75th percentile.}
      \label{fig:spreadRunTimes48}
  \end{figure}

  \begin{figure}[htp]
      \centering
      \includegraphics[trim={1.0cm 0.8cm 2.0cm 1.5cm},clip,width=0.80\linewidth]{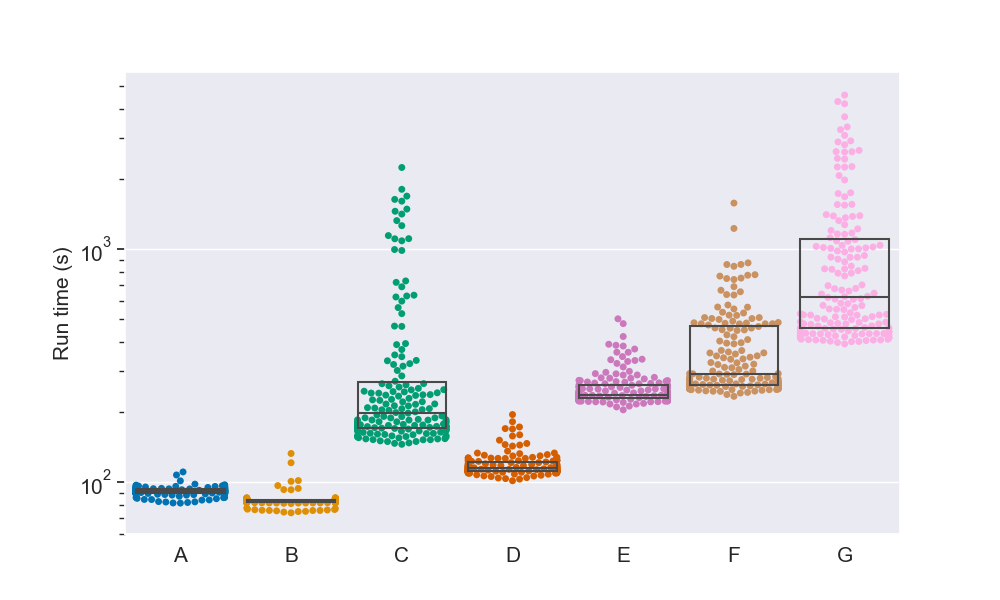}
      \caption{Run times of instances with 96 time steps, sorted by station. One dot represents one instance, however dots may overlap when too many have similar run time. In addition we draw a typical whiskerless box plot. Therefore the lines represent the 25th percentile, the median and the 75th percentile.}
      \label{fig:spreadRunTimes96}
  \end{figure}

  \begin{figure}[htp]
    \centering

    \begin{subfigure}[c]{0.48\textwidth}
    \includegraphics[trim={1.5cm 0.5cm 1.5cm 0.0cm},clip,width=\textwidth]{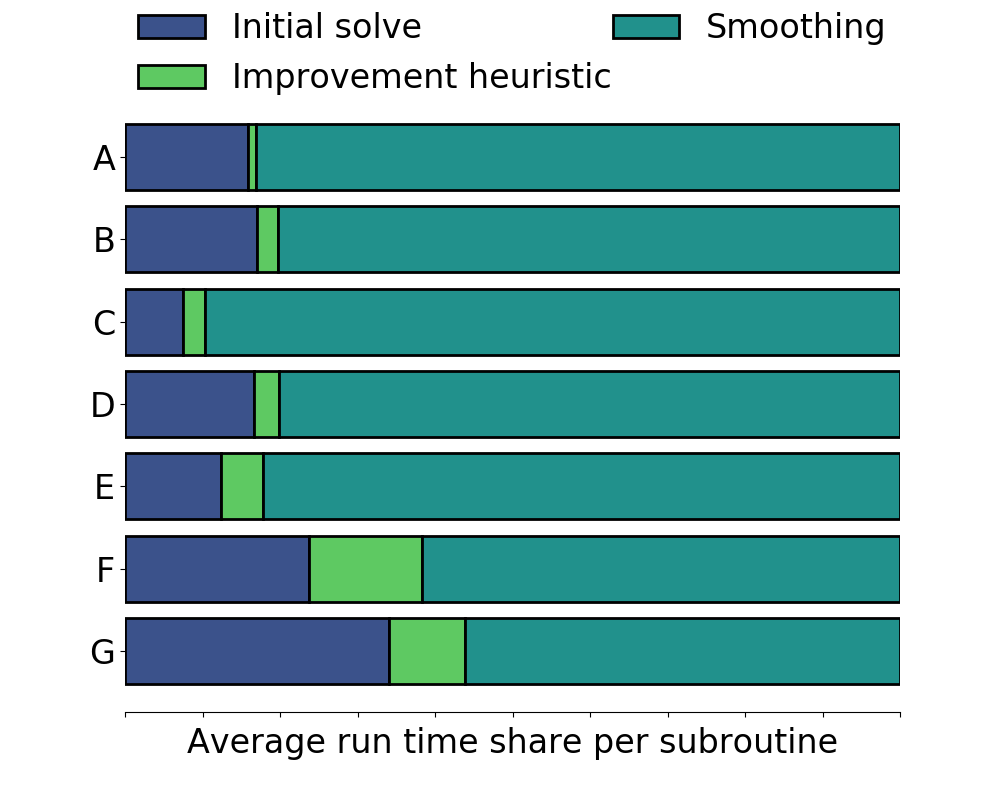}
    \subcaption{24 Timesteps}
    \label{fig:runTimeShare24}
    \end{subfigure}
    \begin{subfigure}[c]{0.48\textwidth}
    \includegraphics[trim={1.5cm 0.5cm 1.5cm 0.0cm},clip,width=\textwidth]{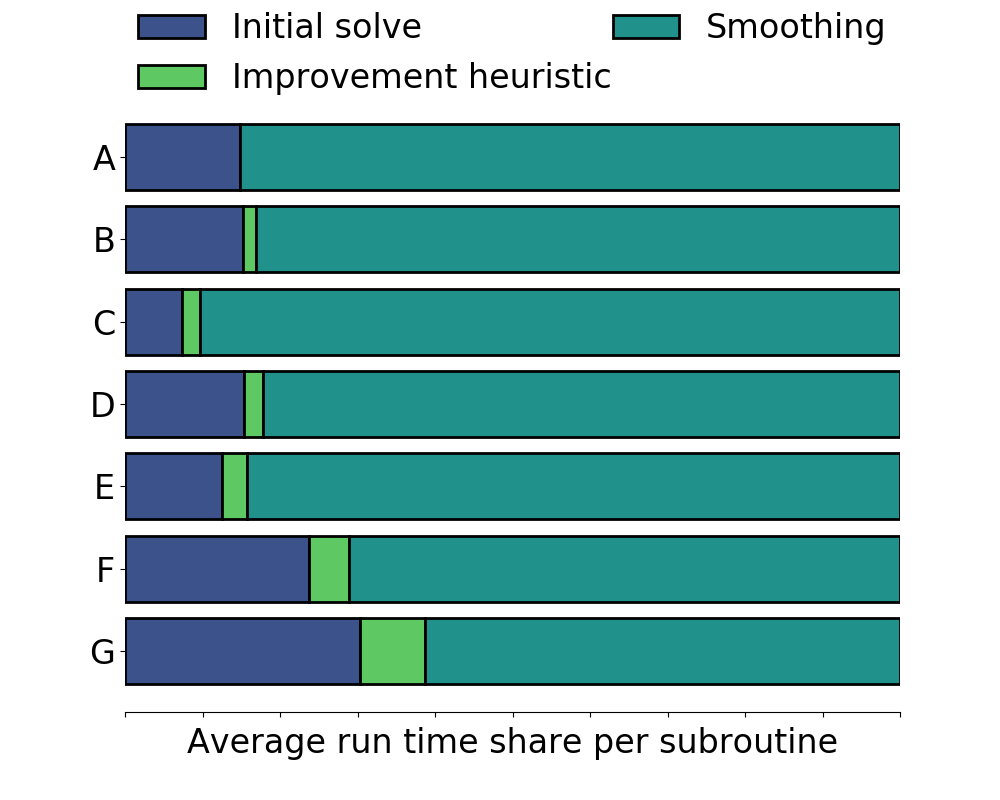}
    \subcaption{48 Timesteps}
    \label{fig:runTimeShare48}
    \end{subfigure}
    \caption{Portion of run time spend in the single subroutines of Algorithm~\ref{algo:completeStationModel} displayed for each of the single stations. Figure~\ref{fig:runTimeShare24} represents the results for 24 time steps, Figure~\ref{fig:runTimeShare48} the results for 48 time steps.}
    \label{fig:runTimeShares2448}
  \end{figure}

\end{document}